\RequirePackage{rotating}

\makeatletter
\def\@cons#1#2{\begingroup\let\@elt\relax\xdef#1{\ifx#1\relax\else#1\fi\@elt #2}\endgroup}
\makeatother
\documentclass{svjour3}

\usepackage{amsmath}
\usepackage{amsfonts}
\usepackage{amssymb}
\usepackage{longtable}
\usepackage{rotating}
\usepackage{multirow}
\usepackage[linesnumbered,ruled,commentsnumbered]{algorithm2e}
\usepackage{algorithmicx}
\usepackage{pdflscape}
\usepackage{booktabs}
\usepackage{hyperref}
\usepackage{cleveref}
\crefname{figure}{fig.}{figs.}%
\usepackage[toc,page]{appendix}

\usepackage{graphicx}
\usepackage{tikz}
\usepackage{pgfplots}
\usepgfplotslibrary{groupplots}
\pgfplotsset{compat=newest}
\usetikzlibrary{arrows}
\usepgfplotslibrary{groupplots}
\usetikzlibrary{decorations}
\usetikzlibrary{decorations.markings}
\usetikzlibrary{shapes}
\usetikzlibrary{intersections}
\usetikzlibrary{patterns}

\usepackage{relsize} 

\makeatletter \let\cl@chapter\relax \makeatother

\makeatletter
\newcommand{\labelx}[1]{
	\relax
	\ifmmode
	\label{#1}
	\else
	\ifnum\pdfstrcmp{\@currenvir}{document}=0
	\label{#1}
	\else
	\label[\@currenvir]{#1}
	\fi
	\fi
}
\makeatother

\pgfplotsset{compat=1.16}

\newcommand{\direct}{\texttt{DIRECT}}
\newcommand{\directgl}{\texttt{DIRECT-GL}}

\newcommand{\directav}{\texttt{DIRECT-a}}
\newcommand{\directm}{\texttt{DIRECT-m}}
\newcommand{\birect}{\texttt{BIRECT}}

\newcommand{\adc}{\texttt{ADC}}

\newcommand{\aggresive}{\texttt{Aggressive DIRECT}}
\newcommand{\indexsett}{\mathbb{I}}
\newcommand{\directlib}{\texttt{DIRECTGOLib v1.1}}
\newcommand{\directgo}{\texttt{DIRECTGO v1.1.0}}

\newcommand{\algrule}[1][.2pt]{\par\vskip.5\baselineskip\hrule height #1\par\vskip.5\baselineskip}

\newcommand{\nosemic}{\renewcommand{\@endalgocfline}{\relax}}
\newcommand{\dosemic}{\renewcommand{\@endalgocfline}{\algocf@endline}}
\let\oldnl\nl
\newcommand{\nonl}{\renewcommand{\nl}{\let\nl\oldnl}}

\definecolor{redwood}{rgb}{0.67, 0.31, 0.32}

\emergencystretch=\maxdimen
\definecolor{onyx}{rgb}{0.06, 0.06, 0.06}
\definecolor{sandstorm}{rgb}{0.93, 0.84, 0.25}
\definecolor{princetonorange}{rgb}{1.0, 0.56, 0.0}
\definecolor{sienna}{rgb}{0.53, 0.18, 0.09}
\definecolor{psychedelicpurple}{rgb}{0.87, 0.0, 1.0}

\begin{document}

	\title{An empirical study of various candidate selection and partitioning techniques in the \direct{} framework
	}

	\titlerunning{Study of various candidate selection and partitioning techniques used in \direct}        

	\author{Linas Stripinis         \and
		Remigijus Paulavi\v{c}ius 
	}


	\institute{L. Stripinis, R. Paulavi\v{c}ius \at
		Vilnius University, Institute of Data Science and Digital Technologies, Akademijos 4, LT-08663 Vilnius, Lithuania \\
		\email{linas.stripinis@mif.vu.lt}           
		\and
		R. Paulavi\v{c}ius \at
		\email{remigijus.paulavicius@mif.vu.lt}
	}

	\date{Received: date / Accepted: date}

	\maketitle

	\begin{abstract}
		Over the last three decades, many attempts have been made to improve the \direct{} (DIviding RECTangles) algorithm's efficiency.
		Various novel ideas and extensions have been suggested.
		The main two steps of \direct-type algorithms are selecting and partitioning potentially optimal rectangles.
		However, the most efficient combination of these two steps is an area that has not been investigated so far.
    This paper presents a study covering an extensive examination of various candidate selection and partitioning techniques within the same \direct{} algorithmic framework.
		Twelve \direct-type algorithmic variations are compared on $800$  randomly generated GKLS-type test problems and $96$ box-constrained global optimization problems from \directlib{} with varying complexity.
		Based on these studies, we have identified the most efficient selection and partitioning combinations leading to new, more efficient, \direct-type algorithms.
    All these algorithms are included in the latest version of \directgo{} and are publicly available.
	\end{abstract}

\section{Introduction}
\label{intro}
In this paper, we consider a box-constrained global optimization problem of the form:
\begin{equation}
	\label{eq:opt-problem}
	\begin{aligned}
		& \min_{\mathbf{x}\in D} && f(\mathbf{x}) 
	\end{aligned}
\end{equation}
where $f:\mathbb{R}^n \rightarrow \mathbb{R}$ is a Lipschitz-continuous, potentially ``black-box''  objective function, and $\mathbf{x}$ is the input vector.
Thus, we assume that the analytical information of the objective function $f$ is unknown and can only be obtained by evaluating $f$ at various points of the feasible region, which is an $n$-dimensional hyper-rectangle
\[ 
D = [ \mathbf{a},  \mathbf{b}] = \{ \mathbf{x} \in \mathbb{R}^n: a_j \leq x_j \leq b_j, j = 1, \dots, n\}.
\]
Moreover, $f$ can be non-linear, multi-modal, non-convex, and non-differentiable.

The optimization community attracted considerable interest from the simplicity and efficiency of the deterministic \direct-type algorithms.
The original \direct{} algorithm was developed by Jones et al.~\cite{Jones1993} and is a well-known and widely used solution technique for derivative-free global optimization.
The \direct{} algorithm extends classical Lipschitz optimization~\cite{Paulavicius2006,Paulavicius2007,Paulavicius2008,Paulavicius2009b,Pinter1996book,Piyavskii1967,Sergeyev2011,Shubert1972}, where the need for the Lipschitz constant is eliminated.
This feature made \direct-type methods especially attractive for solving various real-world optimization problems~(see, e.g., \cite{Baker2000,Bartholomew2002,Carter2001,Cox2001,Serafino2011,Gablonsky2001,Liuzzi2010,Paulavicius2019:eswa,Paulavicius2014:book,Stripinis2018b} and the references given therein).
Furthermore, the extensive numerical benchmarks in~\cite{Rios2013} revealed an encouraging performance of the \direct{} algorithm among other tested derivative-free global optimization approaches, belonging to genetic~\cite{John1975}, simulated annealing~\cite{Kirkpatrick1983}, and particle swarm optimization~\cite{Kennedy1995}.

\begin{figure}[ht]
	\resizebox{\textwidth}{!}{
		\begin{tikzpicture}
			\begin{axis}[
				width=0.45\textwidth,height=0.45\textwidth, 
				xlabel = {$c_1$},
				ylabel = {$c_2$},
				enlargelimits=0.05,
				title={Initialization},
				legend style={draw=none},
				legend columns=1, 
				legend style={at={(0.925,-0.2)},font=\normalsize},
				ylabel style={yshift=-0.1cm},
				xlabel style={yshift=0.1cm},
				ytick distance=1/6,
				xtick distance=1/6,
				every axis/.append style={font=\normalsize},
				yticklabels={$0$, $0$,$\frac{1}{6}$, $\frac{1}{3}$, $\frac{1}{2}$, $\frac{2}{3}$, $\frac{5}{6}$, $1$},
				xticklabels={$0$, $0$,$\frac{1}{6}$, $\frac{1}{3}$, $\frac{1}{2}$, $\frac{2}{3}$, $\frac{5}{6}$, $1$},
				]
				\addlegendimage{only marks,mark=*,color=black}
				\addlegendentry{Sampling point}
				\addplot[thick,patch,mesh,draw,black,patch type=rectangle,line width=0.3mm] coordinates {(0,0) (1,0) (1,1) (0,1)} ;
				\draw [black, thick, mark size=0.05pt, fill=blue!50,opacity=0.4,line width=0.3mm] (axis cs:0,0) rectangle (axis cs:1,1);
				\addplot[only marks,mark=*, mark size=1.5pt,black] coordinates {(0.5,0.5)} node[yshift=-8pt] {\tiny $24.13$};
			\end{axis}
		\end{tikzpicture}
		\begin{tikzpicture}
			\begin{axis}[
				width=0.45\textwidth,height=0.45\textwidth, 
				xlabel = {$c_1$},
				enlargelimits=0.05,
				title={Iteration $1$},
				legend style={draw=none},
				legend columns=1, 
				legend style={at={(0.925,-0.2)},font=\normalsize},
				ylabel style={yshift=-0.1cm},
				xlabel style={yshift=0.1cm},
				ytick distance=1/6,
				xtick distance=1/6,
				every axis/.append style={font=\normalsize},
				yticklabels={$0$, $0$,$\frac{1}{6}$, $\frac{1}{3}$, $\frac{1}{2}$, $\frac{2}{3}$, $\frac{5}{6}$, $1$},
				xticklabels={$0$, $0$,$\frac{1}{6}$, $\frac{1}{3}$, $\frac{1}{2}$, $\frac{2}{3}$, $\frac{5}{6}$, $1$},
				]
				\addlegendimage{area legend,blue!30,fill=blue!50,opacity=0.4}
				\addlegendentry{Selected POH}
				\addplot[thick,patch,mesh,draw,black,patch type=rectangle,line width=0.3mm] coordinates {(0,0) (1,0) (1,1) (0,1)} ;
				\addplot[thick,patch,mesh,draw,black,patch type=rectangle,line width=0.3mm] coordinates {(0,0) (1,0) (1,0.3333) (0,0.3333)};
				\addplot[thick,patch,mesh,draw,black,patch type=rectangle,line width=0.3mm] coordinates {(0,0) (1,0) (1,0.6666) (0,0.6666)};
				\addplot[thick,patch,mesh,draw,black,patch type=rectangle,line width=0.3mm] coordinates {(0,0.3333) (0,0.6666) (0.3333,0.6666) (0.3333,0.3333)};
				\addplot[thick,patch,mesh,draw,black,patch type=rectangle,line width=0.3mm] coordinates {(1,0.3333) (1,0.6666) (0.6666,0.6666) (0.6666,0.3333)};
				\draw [black, thick, mark size=0.1pt, fill=blue!50,opacity=0.4,line width=0.3mm] (axis cs:0,0) rectangle (axis cs:1,0.3333);
				\addplot[only marks,mark=*, mark size=1.5pt,black] coordinates {(0.5,0.5)} node[yshift=-8pt] {\tiny $24.13$};
				\addplot[only marks,mark=*, mark size=1.5pt,black] coordinates {(0.5,0.16666)} node[yshift=-8pt] {\tiny $2.41$};
				\addplot[only marks,mark=*, mark size=1.5pt,black] coordinates {(0.5,0.83333)} node[yshift=-8pt] {\tiny $95.84$};
				\addplot[only marks,mark=*, mark size=1.5pt,black] coordinates {(0.16666,0.5)} node[yshift=-8pt] {\tiny $13.10$};
				\addplot[only marks,mark=*, mark size=1.5pt,black] coordinates {(0.83333,0.5)} node[yshift=-8pt] {\tiny $51.39$};
			\end{axis}
		\end{tikzpicture}
		\begin{tikzpicture}
			\begin{axis}[
				width=0.45\textwidth,height=0.45\textwidth, 
				xlabel = {$c_1$},
				enlargelimits=0.05,
				title={Iteration $2$},
				legend style={draw=none},
				legend columns=1, 
				legend style={at={(0.925,-0.2)},font=\normalsize},
				ylabel style={yshift=-0.1cm},
				xlabel style={yshift=0.1cm},
				ytick distance=1/6,
				xtick distance=1/6,
				every axis/.append style={font=\normalsize},
				yticklabels={$0$, $0$,$\frac{1}{6}$, $\frac{1}{3}$, $\frac{1}{2}$, $\frac{2}{3}$, $\frac{5}{6}$, $1$},
				xticklabels={$0$, $0$,$\frac{1}{6}$, $\frac{1}{3}$, $\frac{1}{2}$, $\frac{2}{3}$, $\frac{5}{6}$, $1$},
				]
				\addlegendimage{area legend,black,fill=white,opacity=0.5}
				\addlegendentry{Unselected region}
				\addplot[thick,patch,mesh,draw,black,patch type=rectangle,line width=0.3mm] coordinates {(0,0) (1,0) (1,1) (0,1)} ;
				\draw [black, thick, mark size=0.1pt,line width=0.3mm] (axis cs:0.3333,0.3333) rectangle (axis cs:0.6666,0.6666);
				\draw [black, thick, mark size=0.1pt,line width=0.3mm] (axis cs:0.3333,0) rectangle (axis cs:0.6666,0.6666);
				\draw [black, thick, mark size=0.1pt, fill=blue!50,opacity=0.4,line width=0.3mm] (axis cs:0,0.6666) rectangle (axis cs:1,1);
				\draw [black, thick, mark size=0.1pt, fill=blue!50,opacity=0.4,line width=0.3mm] (axis cs:0.3333,0) rectangle (axis cs:0.6666,0.3333);
				\draw [black, thick, mark size=0.1pt,line width=0.3mm] (axis cs:0.6666,0.3333) rectangle (axis cs:1,0.6666);
				\draw [black, thick, mark size=0.1pt,line width=0.3mm] (axis cs:0,0.6666) rectangle (axis cs:1,1);
				\draw [black, thick, mark size=0.1pt,line width=0.3mm] (axis cs:0,0.3333) rectangle (axis cs:0.3333,0.6666);
				
				\addplot[only marks,mark=*, mark size=1.5pt,black] coordinates {(0.5,0.5)} node[yshift=-8pt] {\tiny $24.13$};
				\addplot[only marks,mark=*, mark size=1.5pt,black] coordinates {(0.5,0.16666)} node[yshift=-8pt] {\tiny $2.41$};
				\addplot[only marks,mark=*, mark size=1.5pt,black] coordinates {(0.5,0.83333)} node[yshift=-8pt] {\tiny $95.84$};
				\addplot[only marks,mark=*, mark size=1.5pt,black] coordinates {(0.16666,0.5)} node[yshift=-8pt] {\tiny $13.10$};
				\addplot[only marks,mark=*, mark size=1.5pt,black] coordinates {(0.83333,0.5)} node[yshift=-8pt] {\tiny $51.39$};
				\addplot[only marks,mark=*, mark size=1.5pt,black] coordinates {(0.16666,0.16666)} node[yshift=-8pt] {\tiny $70.96$};
				\addplot[only marks,mark=*, mark size=1.5pt,black] coordinates {(0.83333,0.16666)} node[yshift=-8pt] {\tiny $14.69$};
			\end{axis}
	\end{tikzpicture}} 
	\caption{Illustration of selection, sampling, and subdivision (trisection) used in the original \direct{} algorithm~\cite{Jones1993} on a two-dimensional \textit{Branin} test function in the first two iterations}
	\label{fig:divide}
\end{figure}
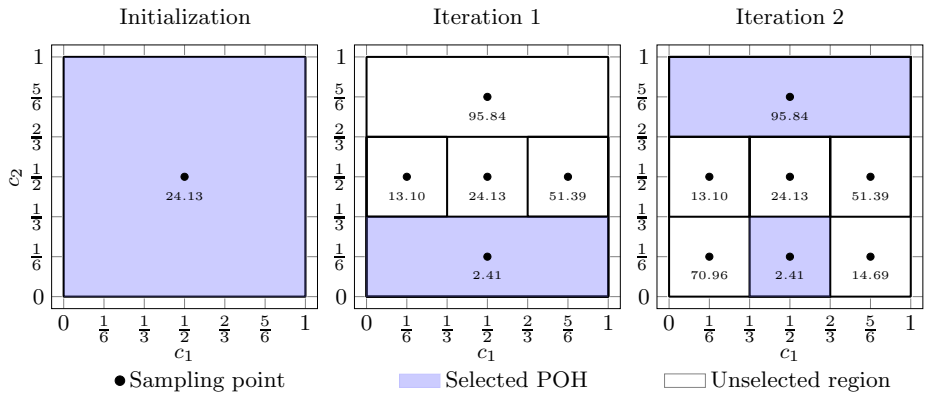

Typically, the \direct-type algorithms include three main steps: selection, sampling, and partitioning (subdivision). 
At each iteration, a specific \direct-type algorithm identifies (selects) the set of potentially optimal hyper-rectangles (POHs) and then samples and subdivides them.
The original \direct{} algorithm uses hyper-rectangular subdivisions based on $n$-dimensional trisection.
The objective function is evaluated at the center points of the newly-formed sub-rectangles. 
Moreover, if several dimensions have the maximum side length, \direct{} starts trisection from the dimension with the lowest $w_j$ and continues to the highest~\cite{Jones2021,Jones1993}.
Here $w_j$ is defined as the best function values sampled along dimension $j$
\begin{equation}
	w_j = \min \{ f(\mathbf{c} + \delta\mathbf{e}_j),  f(\mathbf{c} - \delta\mathbf{e}_j) \},
\end{equation}
where $j \in M$ (set of dimensions with the maximum side length), $\delta$ is equal to one-third of the maximum side length, $\mathbf{c}$ is the center of the hyper-rectangle, and $\mathbf{e}_j$ is the $j$th unit vector.
\Cref{fig:divide} illustrates the selection, sampling, and subdivision (trisection) in the original \direct{} algorithm for a two-dimensional \textit{Branin} test function.

Since the original \direct{} algorithm was published, various \direct-type extensions and modifications have been proposed.
One large group of existing modifications aim to improve the selection of POHs (see, e.g., ~\cite{Baker2000,Gablonsky2001:phd,Mockus2017,Paulavicius2019:eswa,Stripinis2018a}), while the other group concentrates on different partitioning techniques (see, e.g.,~\cite{Jones2001,Liu2015b,Paulavicius2016:jogo,Paulavicius2013:jogo,Sergeyev2006}).
In addition, the authors also make some modifications to the other steps of their algorithms.
Consequently, it is unclear which suggested improvements have the most potential within the \direct{} algorithmic framework.

We address this problem by comparing various proposed candidate selection and partitioning techniques for the remaining algorithmic steps under the same conditions.
This way, we seek to improve the efficiency of existing \direct-type algorithms by creating new combinations based on the previous proposals.
Twelve mostly new  \direct-type algorithmic variations are introduced and investigated using three selection and four partitioning schemes.

The rest of the paper is organized as follows.
\Cref{rewiev} reviews the original \direct{} algorithm and well-known \direct-type modifications proposed for the candidate selection and subdivision.
The obtained new combinations are described in \Cref{sec:combinations}.
An extensive experimental analysis using traditional test problems is presented in \Cref{sec:experiments}, while on GKLS-type test problems in \Cref{sec:exp-GKLS}.
Finally, in \Cref{sec:conclusiuo}, we conclude the paper.

\section{Overview of candidate selection and partitioning techniques used in \direct-type algorithms}
\label{rewiev}

This section reviews the most well-known strategies for selecting and partitioning potentially optimal candidates in \direct-type algorithms.
We start with a brief review of the main steps of the original \direct{} algorithm, with particular emphasis on candidate selection and partitioning techniques.

\subsection{Original \direct{} algorithm}
The original \direct{} algorithm is a deterministic derivative-free global optimization~\cite{Horst1995:book,Sergeyev2017:book,Strongin2000:book} algorithm subject to simple box constraints.
The main steps of \direct{} are summarized in~\Cref{alg:direct}.
At the \textbf{Initialization} step (see~\Cref{alg:direct}, Lines~\ref{alg:initialization_begin}--\ref{alg:initialization_end}), \direct{} normalizes the search region $D$ to unit hyper-rectangle $\bar{D}$ and refers to the original space $D$ only when evaluating the objective function.
Regardless of the dimension $n$, the first evaluation of the objective function is performed at the midpoint of the unit hyper-rectangle $\mathbf{c}^1 = (1/2, \dots,1/2) \in \bar{D}$.

\begin{algorithm}[ht]
	\normalsize
	\LinesNumbered
	\SetAlgoLined
	\SetKwInOut{Input}{input}
	\SetKwInOut{Output}{output}
	
	\SetKwData{Mmax}{M$_{\rm max}$}
	\SetKwData{Kmax}{K$_{\rm max}$}
	\SetKw{And}{and}
	\SetKw{Or}{or}
	
	\direct($f$,$D$,$opt$);\\
	\Input{Objective function $(f)$, search domain $(D)$, and adjustable algorithmic options $(opt)$: tolerance ($\varepsilon_{\rm pe}$), maximal number of function evaluations ($\Mmax$) and algorithmic iterations ($\Kmax$) ; }
	\Output{The best found objective value $(f^{\min})$, solution point $(\mathbf{x}^{\min})$, and record of various performance metrics: percent error $(pe)$, number of iterations $(k)$, number of function evaluations $(m)$;} 
	\algrule
	
	\nonl \textbf{Initialization step:} \\ 
	\textit{Normalize} the search domain $D$ to be the unit hyper-rectangle $\bar{D}$; \label{alg:initialization_begin} \\
	\textit{Evaluate} the objective function at the center point ($\mathbf{c}^1$) of $\bar{D}$ and set: \\
	$\mathbf{c}^1 = \left(\frac{1}{2}, \frac{1}{2}, ..., \frac{1}{2}\right)$; \\
	$x^{\min}_j = \mid b_j - a_j \mid  c^1_j + a_j, j=1, \dots, n$; \tcp*[f]{referring to $D$}\\
	$f^1 = f(\mathbf{x}^{\min})$, $f^{\min} = f^1$; \\
	\textit{Initialize} performance measures: $k=1$, $m=1$, $pe$;\label{alg:initialization_end} \tcp*[f]{\textit{pe} defined in \eqref{eq:pe}}
	
		\While{$pe > \varepsilon_{\rm pe}$ \And $m < \Mmax$ \And $k < \Kmax$ }{
			
			\textbf{Selection step:} \textit{Identify} the set $S_k$ of POHs using \Cref{def:potOptRect} \; \label{alg:selection_begin}
			
			\ForEach{$\bar{D}^j_k \in S_k$}{
				\textbf{Sampling step:} \textit{Evaluate} $f$ at the newly sampled points in $\bar{D}^j_k$; \label{alg:sampling}\\
				\textbf{Subdivision step:} \textit{Trisect} $\bar{D}^j_k$ as illustrated in \Cref{fig:divide} \; \label{alg:subdivision}
			}\label{alg:global_end}
			
			\textit{Update} $f^{\min}, \mathbf{x}^{\min}$, and performance measures: $k$, $m$ and $pe$;
		}
		\textbf{Return} $f^{\min}, \mathbf{x}^{\min}$, and performance measures: $k$, $m$ and $pe$.
		\caption{Main steps of the \direct{} algorithm}
		\label{alg:direct}
	\end{algorithm}
	Two of the most critical steps in the original \direct{} and other existing modifications are \textbf{Selection} and \textbf{Subdivision}.
	
	\subsubsection{Original selection strategy}
	\label{sssec:DIRECT-selection}
	
	Let the current partition at the iteration $ k $ is defined as
	\[
	\mathcal{P}_k = \{ \bar{D}^i_k : i \in \indexsett_k \}, 
	\]
	where $ \bar{D}^i_k = [\mathbf{a}^i, \mathbf{b}^i] = \{ \mathbf{x} \in \bar{D}: 0 \leq a_j^i \leq x_j \leq b_j^i \leq 1, j = 1,\dots, n, \forall i \in \indexsett_k \} $ and $ \indexsett_k $ is the index set identifying the current partition $ \mathcal{P}_k $.
	The next partition, $\mathcal{P}_{k+1}$, is obtained by subdividing selected POHs from the current partition $ \mathcal{P}_k $.
	Note there is only one candidate, $\bar{D}^1_1$, that at the first iteration ($k=1$), which is automatically potentially optimal.
	The formal requirement of potential optimality in subsequent iterations is stated in \Cref{def:potOptRect}.
	
	\begin{definition}
		Let $ \mathbf{c}^i $ denote the center sampling point and $ \delta^i_k $ be a measure (equivalently, sometimes called distance or size) of the hyper-rectangle $ \bar{D}^i_k$.
		Let $ \varepsilon > 0 $ be a positive constant and $f^{\min}$ be the best currently found value of the objective function.
		A hyper-rectangle $ \bar{D}^h_k, h \in \indexsett_k $ is said to be potentially optimal if there exists some rate-of-change (Lipschitz) constant $ \tilde{L} > 0$ such that
		\begin{eqnarray}
			f(\mathbf{x}^h) - \tilde{L}\delta^h_k & \leq & f(\mathbf{x}^i) - \tilde{L}\delta^i_k, \quad \forall i \in \indexsett_k, \label{eqn:potOptRect1} \\
			f(\mathbf{x}^h) - \tilde{L}\delta^h_k & \leq & f^{\min} - \varepsilon|f^{\min}|, \label{eqn:potOptRect2}
		\end{eqnarray}
		where 
		\begin{equation}
			\label{eq:space_original}
			x^i_j = \mid b_j - a_j \mid  c^i_j + a_j, j=1,\dots,n,
		\end{equation} 
		and the measure of the hyper-rectangle $ \bar{D}^i_k$ is
		\begin{equation}
			\label{eq:distance}
			\delta^i_k = \frac{1}{2} \| {\mathbf{b}}^i - {\mathbf{a}}^i \|_2.
		\end{equation}
		\label{def:potOptRect}
	\end{definition}
	
	The hyper-rectangle $ \bar{D}^j_k $ is potentially optimal if the lower Lipschitz bound for the objective function computed by the left-hand side of \eqref{eqn:potOptRect1} is the smallest one with some positive constant $\tilde{L}$ in the current partition $ \mathcal{P}_k $.
	In~\eqref{eqn:potOptRect2}, the parameter $\varepsilon$ is used to protect from an excessive refinement of the local minima~\cite{Jones1993,Paulavicius2014:jogo}.
	In~\cite{Jones1993}, the authors obtained good results for $\varepsilon$ values ranging from $10^{-3}$ to $10^{-7}$.
	A geometrical interpretation of POH selection using \Cref{def:potOptRect} is illustrated on the left panel of \Cref{fig:selection}.
	Here each hyper-rectangle is represented as a point.
	The $x$-axis shows the size of the measure $(\delta)$ while the $y$-axis -- the objective function value attained at the midpoint $(\mathbf{c})$ of this hyper-rectangle.
	The hyper-rectangles meeting conditions \eqref{eqn:potOptRect1} and \eqref{eqn:potOptRect2} are points on the lower-right convex hull (highlighted in blue color).
	
	However, such a selection strategy can be especially inefficient, e.g., for symmetric problems.
	There may be many POHs with the same diameter $\delta^i_k$ and objective value, leading to a drastic increase of selected POHs per iteration.
	To overcome this, authors in \cite{Gablonsky2001} proposed selecting only one of these many ``equivalent'' candidates.
	In \cite{Jones2021}, the authors revealed that such modification could significantly increase the performance of the \direct{} algorithm.
	In this paper, we call this an \textit{improved original selection strategy}.
	
	\subsubsection{Original partitioning scheme}
	\label{sssec:DIRECT-partitioning}
	
	In the \textbf{Sampling} and \textbf{Subdivision} steps (see Algorithm~\ref{alg:direct}, Lines~\ref{alg:sampling} and \ref{alg:subdivision}), a hyper-rectangular partition based on $n$-dimensional trisection is used.
	Using this scheme, the POHs are partitioned into smaller non-intersecting hyper-rectangles (see~\Cref{fig:divide}), containing the lower function values in larger new hyper-rectangles.
	
\subsection{Other candidate selection schemes in \direct-type algorithms}
\label{modifications2}

Various improvements and new ideas for candidate selection were proposed in the literature.
To prevent the \direct{} algorithm from being sensitive to the objective function's additive scaling, authors in~\cite{Finkel2006} introduced a scaling of the objective function values by subtracting the median value calculated from the previously evaluated function values.
More specifically, in the selection step, a new \directm{} replaces \eqref{eqn:potOptRect2} from \Cref{eqn:potOptRect2} to: 
\begin{equation}
	f(\mathbf{x}^i) - \tilde{L}\delta^i \leq f^{\min} - \varepsilon|f^{\min} - {f}^{\rm median}|. \label{eqn:potOptRectz}
\end{equation}
Similarly, in~\cite{Liu2013}, authors adopted the similar idea in \directav{}.
At each iteration, instead of the median value (${f}^{\rm median}$), authors proposed to use the average value $({f}^{\rm average})$:
\begin{equation}
	f(\mathbf{x}^i) - \tilde{L}\delta^i \leq f^{\min} - \varepsilon|f^{\min} - {f}^{\rm average}|. \label{eqn:potOptRect4}
\end{equation}
The authors in~\cite{Finkel2004aa,Liu2015b,Liu2015} showed that different schemes controlling the $\varepsilon$ parameter in \eqref{eqn:potOptRect2} could increase the efficiency of the \direct{} algorithm, especially when needed to fine-tune the solution to higher accuracy.

In order to verify this, the experimental investigation of the original \direct, \directm{} (based on \cref{eqn:potOptRectz}), and \directav{} (based on \cref{eqn:potOptRect4}) algorithms on an extensive set consisting of $81$ test and six engineering problems (from \texttt{DIRECTGOLib v1.0} \cite{DIRECTGOLib2022v10}) were performed in \cite{Stripinis2021:dgo}.
Our investigation revealed that in solving engineering problems, a significant performance difference was not observed.
However, on $81$ test problems, the original \direct{} proved to be more efficient, and using the same stopping conditions solved $10$ and $23$ more test problems than \directm{} and \directav{} accordingly (for details, see Table $3$ in \cite{Stripinis2021:dgo}).
Based on this, the original eq. \eqref{eqn:potOptRect2} was used in an improved original selection strategy in our experimental study.


Below we focus on the other two selection schemes considered in this research.

\subsubsection{Aggressive selection}
\label{sssec:aggressive-selection}

In~\cite{Baker2000}, the authors relaxed the selection criteria of POHs and proposed an aggressive version of the \direct{} algorithm (\aggresive).
The main idea is to select and divide at least one hyper-rectangle from each group of different diameters $(\delta_k^i)$ containing the lowest objective function value.
Such aggressive selection ensures much more objective function evaluations per iteration compared to other existing POH selection schemes.
From the optimization point of view, such an approach may seem less favorable since it ``wastes'' function evaluations by exploring unnecessary (non-potentially optimal) hyper-rectangles.
However, such a strategy is much more appealing in a parallel environment, as was shown in~\cite{He2009part2,He2009part1,He2010,Watson2001}.

\begin{figure}[htb]
	\resizebox{\textwidth}{!}{
		\begin{tikzpicture}
			\begin{groupplot}[
				group style={
					group size=3 by 1,
					x descriptions at=edge bottom,
					y descriptions at=edge left,
					vertical sep=0pt,
					horizontal sep=6pt,
				},
				height=0.5\textwidth,width=0.7\textwidth,
				]
				\nextgroupplot[
				title  = {\parbox{3.7cm}{Original \direct{} selection}},
				xlabel = {\large $\delta$},
				ylabel = {\large Function values},
				ymin=0,
				ymax=1,
				xmin=0,
				xmax=0.6,
				ytick distance=0.2,
				xtick distance=0.2,
				legend style={font=\scriptsize},
				height=0.6\textwidth,width=0.5\textwidth,
				enlargelimits=0.05,
				]
				\addlegendimage{black!60,only marks,mark=*,mark size=2pt}
				\addlegendentry{non-potentially optimal}
				\addlegendimage{blue,mark=*,mark size=2.5pt}
				\addlegendentry{potentially optimal}
				\addplot[black!60,only marks,mark=*,mark size=2pt] table[x=D,y=F] {data/ddr2.txt};
				\addplot[blue,mark=*,mark size=2.5pt] table[x=D_pot_opt,y=F_pot_opt] {data/ddr2.txt};
				\nextgroupplot[
				title  = {\parbox{3.7cm}{Aggressive selection}},
				xlabel = {\large $\delta$},
				ymin=0,
				ymax=1,
				xmin=0,
				xmax=0.6,
				ytick distance=0.25,
				xtick distance=0.2,
				height=0.6\textwidth,width=0.5\textwidth,
				enlargelimits=0.05,
				]
				\addlegendimage{black!60,only marks,mark=*,mark size=2pt}
				\addlegendentry{non-potentially optimal}
				\addlegendimage{blue,mark=*,mark size=2.5pt}
				\addlegendentry{potentially optimal}
				\addplot[black!60,only marks,mark=*,mark size=2pt] table[x=D,y=F] {data/ddr1.txt};
				\addplot[blue,mark=*,mark size=2.5pt] table[x=D_pot_opt,y=F_pot_opt] {data/ddr1.txt};
				\nextgroupplot[
				title  = {\parbox{3.7cm}{Pareto selection}},
				xlabel = {\large $\delta$},
				ymin=0,
				ymax=1,
				xmin=0,
				xmax=0.6,
				xtick distance=0.2,
				height=0.6\textwidth,width=0.5\textwidth,
				enlargelimits=0.05,
				]
				\addlegendimage{black!60,only marks,mark=*,mark size=2pt}
				\addlegendentry{non-potentially optimal}
				\addlegendimage{blue,mark=*,mark size=2.5pt}
				\addlegendentry{potentially optimal}
				\addplot[black!60,only marks,mark=*,mark size=2pt] table[x=D,y=F] {data/ddr.txt};
				\addplot[blue,mark=*,mark size=2.5pt] table[x=D_pot_opt,y=F_pot_opt] {data/ddr.txt};
			\end{groupplot}
	\end{tikzpicture}}
	\caption{Comparison of three different selection schemes (original \direct, aggressive, and Pareto) applied on the same set of points}
	\label{fig:selection}
\end{figure}
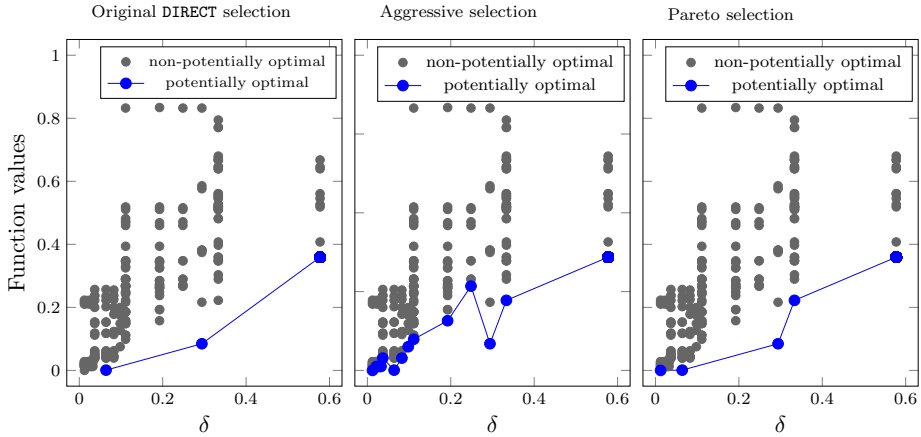

In \cite{He2008}, authors showed that limiting the refinement of the search-space when the size of hyper-rectangles $ \delta^i_k $ reached some prescribed size $ \delta^{\rm limit} $, the memory usage reduces from $10 \%$ to $70 \%$, and the algorithm can run longer without memory allocation failure.
In the experimental part (described in \Cref{sec:experiments}), the limit parameter $ (\delta^{\rm limit}) $ was set to the size of a hyper-rectangle that has been subdivided $ 50 n $ times.
The $ \delta^{\rm limit} $ parameter is intended for the same purpose as the equation \eqref{eqn:potOptRect2} and tries to avoid wasting function estimates by ``over-exploring'' the local minimum area.
We call this an \textit{improved aggressive selection strategy}.
A geometrical interpretation of the aggressive selection is shown in the middle panel of \Cref{fig:selection}.

\subsubsection{Two-step-based Pareto selection}
\label{sssec:two-step-selection}

In a more recent modification \directgl~\cite{Stripinis2018a}, we proposed a new two-step-based selection strategy for the identification of the extended set of POHs. 
In both steps, \directgl{} selects only Pareto optimal hyper-rectangles: in the ﬁrst step, non-dominated on size (the higher, the better) and center point function value (the lower, the better), while in the second, non-dominated on size and distance from the current minimum point (the closer, the better) and takes the unique union of identified candidates in both steps.
We note this scheme does not have any protection against over-exploration in sub-optimal local minima regions.

A geometrical interpretation of the selection procedure is shown in the right panel of \Cref{fig:selection}.
In the first step, \directgl{} selects Pareto hyper-rectangles concerning the size and function value.
Therefore, unlike the \aggresive{} strategy, hyper-rectangles from the groups where the minimum objective function value is higher than the minimum value from the larger groups are not selected in \directgl{}.
Compared to the original selection (Definitions~\ref{def:potOptRect}), in \directgl, the set of POHs is enlarged by adding more medium-sized hyper-rectangles.
In this sense, Pareto selection may be more global than the original \direct{} selection.
Additionally, in the second step, \directgl{} selects the hyper-rectangles that are non-dominated concerning the size and distance from the current minimum point.
This way, the set of POHs is enlarged with various size hyper-rectangles nearest the current minimum point, assuring a broader examination around it.

\subsection{Other partitioning schemes in \direct-type algorithms}
\label{modifications}

\subsubsection{Trisection strategy, along single the longest side}
\label{sssec:trisection-along-single-side}

\begin{figure}[ht]
\resizebox{\textwidth}{!}{
	\begin{tikzpicture}
		\begin{axis}[
			width=0.45\textwidth,height=0.45\textwidth, 
			xlabel = {$c_1$},
			ylabel = {$c_2$},
			enlargelimits=0.05,
			title={Initialization},
			legend style={draw=none},
			legend columns=1, 
			legend style={at={(0.925,-0.2)},font=\normalsize},
			ylabel style={yshift=-0.1cm},
			xlabel style={yshift=0.1cm},
			ytick distance=1/6,
			xtick distance=1/6,
			every axis/.append style={font=\normalsize},
			yticklabels={$0$, $0$,$\frac{1}{6}$, $\frac{1}{3}$, $\frac{1}{2}$, $\frac{2}{3}$, $\frac{5}{6}$, $1$},
			xticklabels={$0$, $0$,$\frac{1}{6}$, $\frac{1}{3}$, $\frac{1}{2}$, $\frac{2}{3}$, $\frac{5}{6}$, $1$},
			]
			\addlegendimage{only marks,mark=*,color=black}
			\addlegendentry{Sampling point}
			\addplot[thick,patch,mesh,draw,black,patch type=rectangle,line width=0.3mm] coordinates {(0,0) (1,0) (1,1) (0,1)} ;
			\draw [black, thick, mark size=0.05pt, fill=blue!50,opacity=0.4,line width=0.3mm] (axis cs:0,0) rectangle (axis cs:1,1);
			\addplot[only marks,mark=*, mark size=1.5pt,black] coordinates {(0.5,0.5)} node[yshift=-8pt] {\tiny $24.13$};
		\end{axis}
	\end{tikzpicture}
	\begin{tikzpicture}
		\begin{axis}[
			width=0.45\textwidth,height=0.45\textwidth, 
			xlabel = {$c_1$},
			enlargelimits=0.05,
			title={Iteration $1$},
			legend style={draw=none},
			legend columns=1, 
			legend style={at={(0.925,-0.2)},font=\normalsize},
			ylabel style={yshift=-0.1cm},
			xlabel style={yshift=0.1cm},
			ytick distance=1/6,
			xtick distance=1/6,
			every axis/.append style={font=\normalsize},
			yticklabels={$0$, $0$,$\frac{1}{6}$, $\frac{1}{3}$, $\frac{1}{2}$, $\frac{2}{3}$, $\frac{5}{6}$, $1$},
			xticklabels={$0$, $0$,$\frac{1}{6}$, $\frac{1}{3}$, $\frac{1}{2}$, $\frac{2}{3}$, $\frac{5}{6}$, $1$},
			]
			\addlegendimage{area legend,blue!30,fill=blue!50,opacity=0.4}
			\addlegendentry{Selected POH}
			\addplot[thick,patch,mesh,draw,black,patch type=rectangle,line width=0.3mm] coordinates {(0,0) (1,0) (1,1) (0,1)} ;
			\draw [black, thick, mark size=0.1pt,line width=0.3mm] (axis cs:0,0) rectangle (axis cs:1,1);
			\draw [black, thick, mark size=0.1pt,line width=0.3mm] (axis cs:0.3333,0) rectangle (axis cs:0.6666,1);
			\draw [black, thick, mark size=0.1pt, fill=blue!50,opacity=0.4,line width=0.3mm] (axis cs:0,0) rectangle (axis cs:0.3333,1);
			\addplot[only marks,mark=*, mark size=1.5pt,black] coordinates {(0.5,0.5)} node[yshift=-8pt] {\tiny $24.13$};
			\addplot[only marks,mark=*, mark size=1.5pt,black] coordinates {(0.16666,0.5)} node[yshift=-8pt] {\tiny $13.10$};
			\addplot[only marks,mark=*, mark size=1.5pt,black] coordinates {(0.83333,0.5)} node[yshift=-8pt] {\tiny $51.39$};
		\end{axis}
	\end{tikzpicture}
	\begin{tikzpicture}
		\begin{axis}[
			width=0.45\textwidth,height=0.45\textwidth, 
			xlabel = {$c_1$},
			enlargelimits=0.05,
			title={Iteration $2$},
			legend style={draw=none},
			legend columns=1, 
			legend style={at={(0.925,-0.2)},font=\normalsize},
			ylabel style={yshift=-0.1cm},
			xlabel style={yshift=0.1cm},
			ytick distance=1/6,
			xtick distance=1/6,
			every axis/.append style={font=\normalsize},
			yticklabels={$0$, $0$,$\frac{1}{6}$, $\frac{1}{3}$, $\frac{1}{2}$, $\frac{2}{3}$, $\frac{5}{6}$, $1$},
			xticklabels={$0$, $0$,$\frac{1}{6}$, $\frac{1}{3}$, $\frac{1}{2}$, $\frac{2}{3}$, $\frac{5}{6}$, $1$},
			]
			\addlegendimage{area legend,black,fill=white,opacity=0.5}
			\addlegendentry{Not selected region}
			\addplot[thick,patch,mesh,draw,black,patch type=rectangle,line width=0.3mm] coordinates {(0,0) (1,0) (1,1) (0,1)} ;
			\draw [black, thick, mark size=0.1pt,line width=0.3mm] (axis cs:0,0) rectangle (axis cs:1,1);
			\draw [black, thick, mark size=0.1pt,line width=0.3mm] (axis cs:0.3333,0) rectangle (axis cs:0.6666,1);
			\draw [black, thick, mark size=0.1pt,line width=0.3mm] (axis cs:0,0) rectangle (axis cs:0.3333,1);
			\draw [black, thick, mark size=0.1pt,line width=0.3mm] (axis cs:0,0.3333) rectangle (axis cs:0.3333,0.6666);
			\draw [black, thick, mark size=0.1pt, fill=blue!50,opacity=0.4,line width=0.3mm] (axis cs:0,0.6666) rectangle (axis cs:0.3333,1);
			\draw [black, thick, mark size=0.1pt, fill=blue!50,opacity=0.4,line width=0.3mm] (axis cs:0.3333,0) rectangle (axis cs:0.6666,1);
			\addplot[only marks,mark=*, mark size=1.5pt,black] coordinates {(0.5,0.5)} node[yshift=-8pt] {\tiny $24.13$};
			\addplot[only marks,mark=*, mark size=1.5pt,black] coordinates {(0.16666,0.5)} node[yshift=-8pt] {\tiny $13.10$};
			\addplot[only marks,mark=*, mark size=1.5pt,black] coordinates {(0.83333,0.5)} node[yshift=-8pt] {\tiny $51.39$};
			\addplot[only marks,mark=*, mark size=1.5pt,black] coordinates {(0.16666,0.16666)} node[yshift=-8pt] {\tiny $70.96$};
			\addplot[only marks,mark=*, mark size=1.5pt,black] coordinates {(0.16666,0.83333)} node[yshift=-8pt] {\tiny $5.24$};
		\end{axis}
\end{tikzpicture}} 
\caption{Two-dimensional illustration of the partitioning technique used in the revised version of the \direct{} algorithm~\cite{Jones2001} on a two-dimensional \textit{Branin} test function}
\label{fig:dividere}
\end{figure}
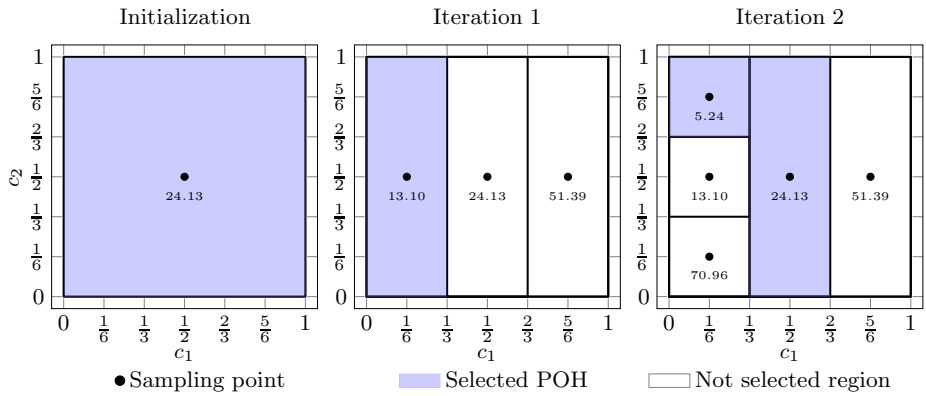

In \cite{Jones2001}, the author proposed a revised version of the original \direct{} algorithm.
One of the main modifications is to trisect selected POHs only along the single longest side (coordinate), see~\Cref{fig:dividere}.
If there are several equal longest sides, the coordinate that has been split the least times during the entire search process so far is selected.
If there is a tie on the latter criterion, the lowest indexed dimension is selected. 
In \cite{Jones2021}, authors showed that dividing a selected rectangle on only one the longest side instead of all can significantly increase the convergence speed.

\subsubsection{Diagonal trisection strategy}
\label{sssec:diagonal-trisection}

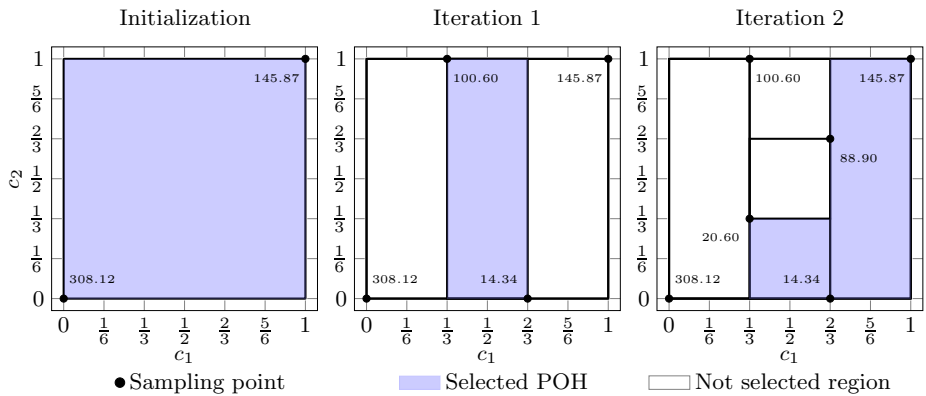
\begin{figure}[ht]
\resizebox{\textwidth}{!}{
	\begin{tikzpicture}
		\begin{axis}[
			width=0.45\textwidth,height=0.45\textwidth, 
			xlabel = {$c_1$},
			ylabel = {$c_2$},
			enlargelimits=0.05,
			title={Initialization},
			legend style={draw=none},
			legend columns=1, 
			legend style={at={(0.925,-0.2)},font=\normalsize},
			ylabel style={yshift=-0.1cm},
			xlabel style={yshift=0.1cm},
			ytick distance=1/6,
			xtick distance=1/6,
			every axis/.append style={font=\normalsize},
			yticklabels={$0$, $0$,$\frac{1}{6}$, $\frac{1}{3}$, $\frac{1}{2}$, $\frac{2}{3}$, $\frac{5}{6}$, $1$},
			xticklabels={$0$, $0$,$\frac{1}{6}$, $\frac{1}{3}$, $\frac{1}{2}$, $\frac{2}{3}$, $\frac{5}{6}$, $1$},
			]
			\addlegendimage{only marks,mark=*,color=black}
			\addlegendentry{Sampling point}
			\addplot[thick,patch,mesh,draw,black,patch type=rectangle,line width=0.3mm] coordinates {(0,0) (1,0) (1,1) (0,1)} ;
			\draw [black, thick, mark size=0.05pt, fill=blue!50,opacity=0.4,line width=0.3mm] (axis cs:0,0) rectangle (axis cs:1,1);
			\addplot[only marks,mark=*, mark size=1.5pt,black] coordinates {(1,1)} node[yshift=-8pt,xshift=-12pt] {\tiny $145.87$};
			\addplot[only marks,mark=*, mark size=1.5pt,black] coordinates {(0,0)} node[yshift=8pt,xshift=12pt] {\tiny $308.12$};
		\end{axis}
	\end{tikzpicture}
	\begin{tikzpicture}
		\begin{axis}[
			width=0.45\textwidth,height=0.45\textwidth, 
			xlabel = {$c_1$},
			enlargelimits=0.05,
			title={Iteration $1$},
			legend style={draw=none},
			legend columns=1, 
			legend style={at={(0.925,-0.2)},font=\normalsize},
			ylabel style={yshift=-0.1cm},
			xlabel style={yshift=0.1cm},
			ytick distance=1/6,
			xtick distance=1/6,
			every axis/.append style={font=\normalsize},
			yticklabels={$0$, $0$,$\frac{1}{6}$, $\frac{1}{3}$, $\frac{1}{2}$, $\frac{2}{3}$, $\frac{5}{6}$, $1$},
			xticklabels={$0$, $0$,$\frac{1}{6}$, $\frac{1}{3}$, $\frac{1}{2}$, $\frac{2}{3}$, $\frac{5}{6}$, $1$},
			]
			\addlegendimage{area legend,blue!30,fill=blue!50,opacity=0.4}
			\addlegendentry{Selected POH}
			\addplot[thick,patch,mesh,draw,black,patch type=rectangle,line width=0.3mm] coordinates {(0,0) (1,0) (1,1) (0,1)} ;
			\draw [black, thick, mark size=0.1pt,line width=0.3mm] (axis cs:0,0) rectangle (axis cs:1,1);
			\draw [black, thick, mark size=0.1pt,line width=0.3mm] (axis cs:0.3333,0) rectangle (axis cs:0.6666,1);
			\draw [black, thick, mark size=0.1pt, fill=blue!50,opacity=0.4,line width=0.3mm] (axis cs:0.3333,0) rectangle (axis cs:0.6666,1);
			\addplot[only marks,mark=*, mark size=1.5pt,black] coordinates {(1,1)} node[yshift=-8pt,xshift=-12pt] {\tiny $145.87$};
			\addplot[only marks,mark=*, mark size=1.5pt,black] coordinates {(0,0)} node[yshift=8pt,xshift=12pt] {\tiny $308.12$};
			\addplot[only marks,mark=*, mark size=1.5pt,black] coordinates {(0.6666,0)} node[yshift=8pt,xshift=-12pt] {\tiny $14.34$};
			\addplot[only marks,mark=*, mark size=1.5pt,black] coordinates {(0.3333,1)} node[yshift=-8pt,xshift=12pt] {\tiny $100.60$};
		\end{axis}
	\end{tikzpicture}
	\begin{tikzpicture}
		\begin{axis}[
			width=0.45\textwidth,height=0.45\textwidth, 
			xlabel = {$c_1$},
			enlargelimits=0.05,
			title={Iteration $2$},
			legend style={draw=none},
			legend columns=1, 
			legend style={at={(0.925,-0.2)},font=\normalsize},
			ylabel style={yshift=-0.1cm},
			xlabel style={yshift=0.1cm},
			ytick distance=1/6,
			xtick distance=1/6,
			every axis/.append style={font=\normalsize},
			yticklabels={$0$, $0$,$\frac{1}{6}$, $\frac{1}{3}$, $\frac{1}{2}$, $\frac{2}{3}$, $\frac{5}{6}$, $1$},
			xticklabels={$0$, $0$,$\frac{1}{6}$, $\frac{1}{3}$, $\frac{1}{2}$, $\frac{2}{3}$, $\frac{5}{6}$, $1$},
			]
			\addlegendimage{area legend,black,fill=white,opacity=0.5}
			\addlegendentry{Not selected region}
			\addplot[thick,patch,mesh,draw,black,patch type=rectangle,line width=0.3mm] coordinates {(0,0) (1,0) (1,1) (0,1)} ;
			\draw [black, thick, mark size=0.1pt,line width=0.3mm] (axis cs:0,0) rectangle (axis cs:1,1);
			\draw [black, thick, mark size=0.1pt,line width=0.3mm] (axis cs:0.3333,0) rectangle (axis cs:0.6666,1);
			\draw [black, thick, mark size=0.1pt,line width=0.3mm] (axis cs:0,0) rectangle (axis cs:0.3333,1);
			\draw [black, thick, mark size=0.1pt,line width=0.3mm] (axis cs:0.3333,0.3333) rectangle (axis cs:0.6666,0.6666);
			\draw [black, thick, mark size=0.1pt, fill=blue!50,opacity=0.4,line width=0.3mm] (axis cs:0.3333,0) rectangle (axis cs:0.6666,0.3333);
			\draw [black, thick, mark size=0.1pt, fill=blue!50,opacity=0.4,line width=0.3mm] (axis cs:0.6666,0) rectangle (axis cs:1,1);
			
			\addplot[only marks,mark=*, mark size=1.5pt,black] coordinates {(1,1)} node[yshift=-8pt,xshift=-12pt] {\tiny $145.87$};
			\addplot[only marks,mark=*, mark size=1.5pt,black] coordinates {(0,0)} node[yshift=8pt,xshift=12pt] {\tiny $308.12$};
			\addplot[only marks,mark=*, mark size=1.5pt,black] coordinates {(0.6666,0)} node[yshift=8pt,xshift=-12pt] {\tiny $14.34$};
			\addplot[only marks,mark=*, mark size=1.5pt,black] coordinates {(0.3333,1)} node[yshift=-8pt,xshift=12pt] {\tiny $100.60$};
			\addplot[only marks,mark=*, mark size=1.5pt,black] coordinates {(0.6666,0.6666)} node[yshift=-8pt,xshift=12pt] {\tiny $88.90$};
			\addplot[only marks,mark=*, mark size=1.5pt,black] coordinates {(0.3333,0.3333)} node[yshift=-8pt,xshift=-12pt] {\tiny $20.60$};
		\end{axis}
\end{tikzpicture}}
\caption{Two-dimensional illustration of the diagonal trisection strategy introduced in the \adc~\cite{Sergeyev2006} algorithm on a two-dimensional \textit{Branin} test function}
\label{fig:divide_adc}
\end{figure}

Adaptive diagonal curves (\adc) based algorithm was introduced in \cite{Sergeyev2006}. 
Independently of the problem dimension, the \adc{} algorithm evaluates the objective function $f(\mathbf{x})$ at two vertices of the main diagonals of each hyper-rectangle $ \bar{D}_k^i$, as shown in \Cref{fig:divide_adc}.
Same as in the revised version of \direct{} \cite{Jones2001}, each selected POH is trisected along just one of the longest sides.
Such a diagonal scheme potentially obtains more comprehensive information about the objective function than center sampling.
The center sampling strategies may sometimes take many iterations to find the solution when a hyper-rectangle containing the optimum has a midpoint with a very bad function value, which makes it undesirable for further selection.
The \adc{} algorithm intuitively reduces this chance for both sampling points in the hyper-rectangle containing the optimum solution by sampling two points per hyper-rectangle.
Therefore, better performance could be expected, especially solving more complex problems.

The main advantage of such a strategy is that it addresses one of the well-known algorithmic weaknesses of the original \direct. 
The feasible region boundary points can only be approached arbitrarily closely but never sampled using the center sampling technique.
Authors in \cite{Huyer1999,Liu2015b} have shown that the latter fact can cause very slow convergence to an optimum if it lies on the feasible region's boundary.

\subsubsection{Diagonal bisection strategy}
\label{sssec:diagonal-bisection}

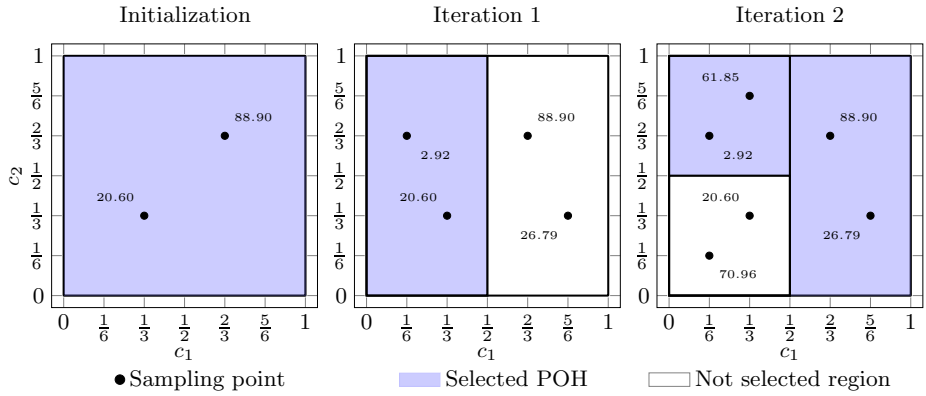
\begin{figure}[ht]
\resizebox{\textwidth}{!}{
	\begin{tikzpicture}
		\begin{axis}[
			width=0.45\textwidth,height=0.45\textwidth, 
			xlabel = {$c_1$},
			ylabel = {$c_2$},
			enlargelimits=0.05,
			title={Initialization},
			legend style={draw=none},
			legend columns=1, 
			legend style={at={(0.925,-0.2)},font=\normalsize},
			ylabel style={yshift=-0.1cm},
			xlabel style={yshift=0.1cm},
			ytick distance=1/6,
			xtick distance=1/6,
			every axis/.append style={font=\normalsize},
			yticklabels={$0$, $0$,$\frac{1}{6}$, $\frac{1}{3}$, $\frac{1}{2}$, $\frac{2}{3}$, $\frac{5}{6}$, $1$},
			xticklabels={$0$, $0$,$\frac{1}{6}$, $\frac{1}{3}$, $\frac{1}{2}$, $\frac{2}{3}$, $\frac{5}{6}$, $1$},
			]
			\addlegendimage{only marks,mark=*,color=black}
			\addlegendentry{Sampling point}
			\addplot[thick,patch,mesh,draw,black,patch type=rectangle,line width=0.3mm] coordinates {(0,0) (1,0) (1,1) (0,1)} ;
			\draw [black, thick, mark size=0.05pt, fill=blue!50,opacity=0.4,line width=0.3mm] (axis cs:0,0) rectangle (axis cs:1,1);
			\addplot[only marks,mark=*, mark size=1.5pt,black] coordinates {(0.3333,0.3333)} node[yshift=8pt,xshift=-12pt] {\tiny $20.60$};
			\addplot[only marks,mark=*, mark size=1.5pt,black] coordinates {(0.6666,0.6666)} node[yshift=8pt,xshift=12pt] {\tiny $88.90$};
		\end{axis}
	\end{tikzpicture}
	\begin{tikzpicture}
		\begin{axis}[
			width=0.45\textwidth,height=0.45\textwidth, 
			xlabel = {$c_1$},
			enlargelimits=0.05,
			title={Iteration $1$},
			legend style={draw=none},
			legend columns=1, 
			legend style={at={(0.925,-0.2)},font=\normalsize},
			ylabel style={yshift=-0.1cm},
			xlabel style={yshift=0.1cm},
			ytick distance=1/6,
			xtick distance=1/6,
			every axis/.append style={font=\normalsize},
			yticklabels={$0$, $0$,$\frac{1}{6}$, $\frac{1}{3}$, $\frac{1}{2}$, $\frac{2}{3}$, $\frac{5}{6}$, $1$},
			xticklabels={$0$, $0$,$\frac{1}{6}$, $\frac{1}{3}$, $\frac{1}{2}$, $\frac{2}{3}$, $\frac{5}{6}$, $1$},
			]
			\addlegendimage{area legend,blue!30,fill=blue!50,opacity=0.4}
			\addlegendentry{Selected POH}
			\addplot[thick,patch,mesh,draw,black,patch type=rectangle,line width=0.3mm] coordinates {(0,0) (1,0) (1,1) (0,1)} ;
			\draw [black, thick, mark size=0.05pt, fill=blue!50,opacity=0.4,line width=0.3mm] (axis cs:0,0) rectangle (axis cs:0.5,1);
			\draw [black, thick, mark size=0.1pt,line width=0.3mm] (axis cs:0,0) rectangle (axis cs:0.5,1);
			\addplot[only marks,mark=*, mark size=1.5pt,black] coordinates {(0.3333,0.3333)} node[yshift=8pt,xshift=-12pt] {\tiny $20.60$};
			\addplot[only marks,mark=*, mark size=1.5pt,black] coordinates {(0.6666,0.6666)} node[yshift=8pt,xshift=12pt] {\tiny $88.90$};
			\addplot[only marks,mark=*, mark size=1.5pt,black] coordinates {(0.8333,0.3333)} node[yshift=-8pt,xshift=-12pt] {\tiny $26.79$};
			\addplot[only marks,mark=*, mark size=1.5pt,black] coordinates {(0.1666,0.6666)} node[yshift=-8pt,xshift=12pt] {\tiny $2.92$};
		\end{axis}
	\end{tikzpicture}
	\begin{tikzpicture}
		\begin{axis}[
			width=0.45\textwidth,height=0.45\textwidth, 
			xlabel = {$c_1$},
			enlargelimits=0.05,
			title={Iteration $2$},
			legend style={draw=none},
			legend columns=1, 
			legend style={at={(0.925,-0.2)},font=\normalsize},
			ylabel style={yshift=-0.1cm},
			xlabel style={yshift=0.1cm},
			ytick distance=1/6,
			xtick distance=1/6,
			every axis/.append style={font=\normalsize},
			yticklabels={$0$, $0$,$\frac{1}{6}$, $\frac{1}{3}$, $\frac{1}{2}$, $\frac{2}{3}$, $\frac{5}{6}$, $1$},
			xticklabels={$0$, $0$,$\frac{1}{6}$, $\frac{1}{3}$, $\frac{1}{2}$, $\frac{2}{3}$, $\frac{5}{6}$, $1$},
			]
			\addlegendimage{area legend,black,fill=white,opacity=0.5}
			\addlegendentry{Not selected region}
			\addplot[thick,patch,mesh,draw,black,patch type=rectangle,line width=0.3mm] coordinates {(0,0) (1,0) (1,1) (0,1)} ;
			\draw [black, thick, mark size=0.1pt, fill=blue!50,opacity=0.4,line width=0.3mm] (axis cs:0.5,0) rectangle (axis cs:1,1);
			\draw [black, thick, mark size=0.1pt, fill=blue!50,opacity=0.4,line width=0.3mm] (axis cs:0,0.5) rectangle (axis cs:0.5,1);
			\draw [black, thick, mark size=0.1pt ,line width=0.3mm] (axis cs:0,0) rectangle (axis cs:0.5,1);
			\draw [black, thick, mark size=0.1pt ,line width=0.3mm] (axis cs:0,0) rectangle (axis cs:0.5,0.5);
			\addplot[only marks,mark=*, mark size=1.5pt,black] coordinates {(0.3333,0.3333)} node[yshift=8pt,xshift=-12pt] {\tiny $20.60$};
			\addplot[only marks,mark=*, mark size=1.5pt,black] coordinates {(0.6666,0.6666)} node[yshift=8pt,xshift=12pt] {\tiny $88.90$};
			\addplot[only marks,mark=*, mark size=1.5pt,black] coordinates {(0.8333,0.3333)} node[yshift=-8pt,xshift=-12pt] {\tiny $26.79$};
			\addplot[only marks,mark=*, mark size=1.5pt,black] coordinates {(0.1666,0.6666)} node[yshift=-8pt,xshift=12pt] {\tiny $2.92$};
			
			\addplot[only marks,mark=*, mark size=1.5pt,black] coordinates {(0.1666,0.1666)} node[yshift=-8pt,xshift=12pt] {\tiny $70.96$};
			\addplot[only marks,mark=*, mark size=1.5pt,black] coordinates {(0.3333,0.8333)} node[yshift=8pt,xshift=-12pt] {\tiny $61.85$};
		\end{axis}
\end{tikzpicture}}
\caption{Two-dimensional illustration of the diagonal bisection strategy used in the \birect{} algorithm~\cite{Paulavicius2016:jogo} on a two-dimensional \textit{Branin} test function}
\label{fig:divide_birect}
\end{figure}

\birect{} (\texttt{BI}secting \texttt{RECT}angles)~\cite{Paulavicius2016:jogo} is motivated by the diagonal partitioning approach~\cite{Sergeyev2006,Sergeyev2008:book,Sergeyev2017:book}.
The bisection is used instead of a trisection typical for diagonal-based and most \direct-type algorithms.
However, neither sampling at the center nor the diagonal's endpoints are appropriate for bisection.
Therefore, in \birect, the objective function is evaluated at two points lying on the diagonal equidistant between themselves and a diagonal's vertices (see \Cref{fig:divide_birect}).
Such a sampling strategy enables the reuse of the sampling points in descendant hyper-rectangles.
Like the \adc{} algorithm (see \Cref{sssec:diagonal-trisection}), \birect{} samples two points per hyper-rectangle.
Therefore, more comprehensive information about the objective function is considered compared to the central sampling strategy used in most \direct-type algorithms.

\section{Summary of new \direct-type algorithmic variations}
\label{sec:combinations}

In this section, we define new variations of the \direct-type algorithms.
In total, twelve variants of \direct-type algorithms are constructed (see \Cref{tab:combo}) by combining three different selection and four partitioning techniques reviewed in the previous section.
\begin{table}[ht]
	\centering
	\caption{Abbreviations of new twelve \direct-type algorithmic variations based on three different selection and four partitioning strategies}
	\label{tab:combo}
	\resizebox{\textwidth}{!}{
		\begin{tabular}{r|r|p{2cm}p{2cm}p{2cm}p{2cm}}
			\toprule
			\multicolumn{2}{c}{} & \multicolumn{4}{c}{\textbf{Partitioning strategy}}\\
			\cmidrule{3-6}
			\multicolumn{2}{c}{} & N-DTC & 1-DTC & 1-DTDV & 1-DBDP \\
			\midrule
			\multirow{3}{*}{\rotatebox{90}{\parbox{1.6cm}{\textbf{Selection scheme}}}} & IO & N-DTC-IO & 1-DTC-IO & 1-DTDV-IO & 1-DBDP-IO \\
			&&&&&\\
			& IA & N-DTC-IA & 1-DTC-IA & 1-DTDV-IA & 1-DBDP-IA \\
			&&&&&\\
			& GL & N-DTC-GL & 1-DTC-GL & 1-DTDV-GL & 1-DBDP-GL \\
			\bottomrule
	\end{tabular}}
\end{table}

The selection strategies used to create these new \direct-type algorithmic variations are:

\begin{itemize}
	\item[1.] \textit{\textbf{I}mproved \textbf{O}riginal selection (\textbf{IO}}) as described in~\Cref{sssec:DIRECT-selection}.
	\item[2.] \textit{\textbf{I}mproved \textbf{A}ggressive selection (\textbf{IA})} as described in \Cref{sssec:aggressive-selection}.
	\item[3.] \textit{Two-step-based (\textbf{G}lobal-\textbf{L}ocal) Pareto selection (\textbf{GL})} as described in~\Cref{sssec:two-step-selection}.
\end{itemize}

The partitioning strategies used in these combinations are the following:
\begin{itemize}
	\item[1.] \textit{Hyper-rectangular partitioning based on \textbf{N}-\textbf{D}imensional \textbf{T}risection and objective function evaluations at \textbf{C}enter points (\textbf{N-DTC})} as described in \Cref{sssec:DIRECT-partitioning}.
	\item[2.] \textit{Hyper-rectangular partitioning based on \textbf{1}-\textbf{D}imensional \textbf{T}risection and objective function evaluations at \textbf{C}enter points (\textbf{1-DTC})} as described in \Cref{sssec:trisection-along-single-side}.
	\item[3.] \textit{Hyper-rectangular partitioning based on \textbf{1}-\textbf{D}imensional \textbf{T}risection and objective function evaluations at two \textbf{D}iagonal \textbf{V}ertices (\textbf{1-DTDV})} as described in \Cref{sssec:diagonal-trisection}.
	\item[4.] \textit{Hyper-rectangular partitioning based on \textbf{1}-\textbf{D}imensional \textbf{B}isection and objective function evaluations at two \textbf{D}iagonals \textbf{P}oints (\textbf{1-DBDP})} as described in \Cref{sssec:diagonal-bisection}.
\end{itemize}

Let us note that some constructed combinations are already used in existing \direct-type algorithms.
For example, the N-DTC-GL algorithm is identical to the recently proposed \directgl{} \cite{Stripinis2018a}.
Furthermore, the N-DTC-IO and 1-DBDP-IO algorithms are highly related to \direct~\cite{Jones1993} and \birect~\cite{Paulavicius2016:jogo} algorithms.
The only difference is that \direct{} and \birect{} algorithms select all ``equivalent'' candidates (with the same diameter and objective function value), while the IO selection rule restricts to one candidate.
The \aggresive~\cite{Baker2000} algorithm is close to the N-DTC-IA variation.
The only difference is the selection step using the limit parameter $ (\delta^{\rm limit}) $.
Moreover, discarding the local search subroutine from the revised hybrid version of the \direct{} algorithm~\cite{Jones2001} would lead to 1-DTC-IO variation.
Finally, the 1-DTDV-IO combination is highly related to the \adc{} algorithm~\cite{Sergeyev2006}, but the latter approach has distinct ``local'' and ``global'' phases in the selection procedure.


\section{Experimental investigation using test problems from \directlib}
\label{sec:experiments}

Test problems from the \directlib{} library~\cite{DIRECTGOLib2022v11} (listed in \Cref{apendixas} \Cref{tab:test}) are used to evaluate the developed algorithms.
In total, we examined new algorithms on $96$  box-constrained global optimization test instances.
Note that different subsets (e.g., low dimensional problems $(n \le 4)$, non-convex problems, etc.) of the entire set were used to deepen the investigation.
All problems and algorithms are implemented in the Matlab R2022a environment and are included in the most recent version of  \direct-type Matlab toolbox \directgo \cite{DIRECTGOv1.1.0}.
All computations were performed on 8th Generation Intel R Core$^\textit{TM}$ i7-8750H @ 2.20GHz Processor.
All $12$ algorithms were tested using a limit of M$_{\rm max} = 10^6$ function evaluations in each run.
For the $96$ analytical test cases with a priori known global optima $ f^* $, the used stopping criterion is based on the percent error:
\begin{equation}
	\label{eq:pe}
	\ pe = 100 \% \times
	\begin{cases}
		\frac{f({\mathbf{x}}) - f^*}{|f^*|},  & f^* \neq 0, \\
		f({\mathbf{x}}),  & f^* = 0,
	\end{cases}
\end{equation}
where $ f^* $ is the known global optimum.
In all experimental studies presented in this section, the algorithms were stopped when the percent error became smaller than the prescribed value $\varepsilon_{\rm pe} = 10^{-2}$ or when the number of function evaluations exceeded the prescribed limit of $10^6$.
In other words, we stop the search when the algorithm has attained an objective function value very close to the known optimum value.

Experimental results presented in this paper are also available in digital form in the \texttt{Results/JOGO} directory of the Github repository~\cite{DIRECTGOv1.1.0}.
The \texttt{Scripts/JOGO} directory of the same Github repository~\cite{DIRECTGOv1.1.0} provides the \texttt{MATLAB} script for cycling through all different classes of \directlib{} test problems used in this paper.
The constructed script can be handy for reproducing the results presented here and comparing and evaluating newly developed algorithms.

\subsection{Investigation of different partitioning strategies}
\label{ssec:partitioning-impact}

First, we compare the performance of new \direct-type algorithms by stressing the used partitioning strategy.
\Cref{tab:aggresive} summarizes comparative results using all twelve \direct-type variations on the whole set of $96$ \directlib{} test problems.
In \Cref{tab:aggresive}, each column corresponds to a different partitioning method.
Since each partitioning method was run on the $96$ problems using $3$ different selection methods (rows of \Cref{tab:aggresive}), it follows that each partitioning method was involved in solving $3 \times 96 = 288$ problems.
The best results are marked in bold.
In total, all three 1-DBDP partitioning-strategy-based algorithms failed to solve ($28/288$) test cases.
In contrast, the second and third best partitioning technique (N-DTC, 1-DTC) based \direct-type approaches did not solve ($29/288$) and ($36/288$) cases accordingly.
Not surprisingly, a higher number of solved test problems leads to a better overall average performance of 1-DBDP partitioning strategy-based \direct-type approaches.
In total, the 1-DBDP partitioning technique-based \direct-type methods required approximately $7 \%$ and $18 \%$ fever function evaluations than the second and third best partitioning scheme (N-DTC and 1-DTC) based \direct-type algorithms accordingly.

However, the situation is different when comparing algorithms based on the median number of function evaluations.
The diagonal trisection strategy (1-DTDV) based \direct-type algorithms are more efficient than competitors.
The median value for all $288$ test problems solved with the 1-DTDV partitioning scheme is approximately $23 \%$ and $47 \%$ better than the second and third best 1-DTC and 1-DBDP approaches accordingly.
Therefore, 1-DTDV partitioning strategy-based \direct-type algorithms can solve at least half of these test problems with the best performance.
Furthermore, most of the time, the 1-DTDV partitioning strategy-based algorithms delivered the best average results in solving low-dimensional test problems ($n \leq 4$).
In total, on $153$ low-dimensional ($n \leq 4$) test cases, the 1-DTDV partitioning strategy-based algorithms required approximately $22 \%$ fever function evaluations than the second-best partition strategy (1-DTC) based variants.
To sum up, the 1-DTDV partitioning strategy performs the best combined with a two-step-based selection scheme (GL).
The 1-DTDV-GL algorithm in general significantly outperformed the other two variations based on IA and IO selection schemes.

The 1-DBDP partitioning strategy-based variation combined with two of the selection schemes (GL and IA) delivered the best average results on higher dimension ($n > 4$) test problems.
In total, the 1-DBDP partitioning strategy-based algorithms required approximately $16 \%$ fever function evaluations than the second-best partition strategy (N-DTC) based variants.

\begin{table}
	\caption{The number of function evaluations of twelve \direct-type variants on \directlib{} test problems}
	\resizebox{\textwidth}{!}{
		\begin{tabular}[tb]{@{\extracolsep{\fill}}lcrrrrrrrr}
			\toprule
			Criteria / Algorithms & $\#$ of cases & \multicolumn{2}{c}{\scriptsize N-DTC-IA} & \multicolumn{2}{c}{\scriptsize 1-DTC-IA} & \multicolumn{2}{c}{\scriptsize 1-DBDP-IA} & \multicolumn{2}{c}{\scriptsize 1-DTDV-IA}  \\
			\midrule
			$\#$ of failed problems & $96$ && $13$ 		&& $13$ 		&& $\mathbf{11}$ 		&& $18$  \\
			Average results			& $96$ && $172,805$ && $160,691$ 	&& $\mathbf{146,887}$ 	&& $202,694$  \\
			Average ($n \leq 4$) 	& $51$ && $25,968$ 	&& $23,638$ 	&& $45,643$ 			&& $\mathbf{9,785}$  \\
			Average ($n > 4$)		& $45$ && $339,791$ && $316,541$ 	&& $\mathbf{262,640}$ 	&& $421,539$  \\
			Average (convex) 		& $30$ && $149,711$ && $126,030$ 	&& $\mathbf{109,374}$ 	&& $153,594$  \\
			Average (non-convex) 	& $66$ && $183,302$ && $176,446$ 	&& $\mathbf{163,939}$ 	&& $225,012$  \\
			Average (uni-modal)		& $15$ && $108,068$ && $78,226$ 	&& $\mathbf{73,957}$ 	&& $111,805$  \\
			Average (multi-modal)	& $81$ && $187,744$ && $179,722$ 	&& $\mathbf{163,717}$ 	&& $223,668$  \\
			Median results			& $96$ && $7,608$ 	&& $\mathbf{1,287}$ && $2,108$ && $1,586$ \\
			\midrule
			Criteria / Algorithms & $\#$ of cases & \multicolumn{2}{c}{\scriptsize N-DTC-IO} & \multicolumn{2}{c}{\scriptsize 1-DTC-IO} & \multicolumn{2}{c}{\scriptsize 1-DBDP-IO} & \multicolumn{2}{c}{\scriptsize 1-DTDV-IO}  \\
			\midrule
			$\#$ of failed problems	& $96$ && $\mathbf{12}$ 		&& $18$ 		&& $\mathbf{12}$ 		&& $21$ 			\\
			Average results			& $96$ && $\mathbf{142,277}$ 	&& $211,463$ 	&& $146,133$ 			&& $227,455$  		\\
			Average ($n \leq 4$) 	& $51$ && $43,832$ 				&& $42,633$ 	&& $\mathbf{41,602}$  	&& $41,990$ \\
			Average ($n > 4$)		& $45$ && $\mathbf{254,819}$ 	&& $403,749$ 	&& $265,522$ 			&& $438,574$  		\\
			Average (convex) 		& $30$ && $111,817$ 			&& $170,675$ 	&& $\mathbf{80,490}$ 	&& $171,868$  		\\
			Average (non-convex) 	& $66$ && $\mathbf{156,122}$ 	&& $230,004$ 	&& $175,971$ 			&& $252,722$  		\\
			Average (uni-modal)		& $15$ && $60,100$ 				&& $\mathbf{57,360}$ && $62,016$ 		&& $111,547$ 		\\
			Average (multi-modal)	& $81$ && $\mathbf{161,240}$ 	&& $247,026$ 	&& $165,545$ 			&& $254,203$  		\\
			Median results			& $96$ && $\mathbf{771}$ 		&& $1,198$ 		&& $953$ 				&& $847$  			\\
			\midrule
			Criteria / Algorithms & $\#$ of cases & \multicolumn{2}{c}{\scriptsize N-DTC-GL} & \multicolumn{2}{c}{\scriptsize 1-DTC-GL} & \multicolumn{2}{c}{\scriptsize 1-DBDP-GL} & \multicolumn{2}{c}{\scriptsize 1-DTDV-GL}  \\
			\midrule
			$\#$ of failed problems & $96$ && $\mathbf{4}$	&& $5$ 				&& $5$ 					&& $5$ \\
			Average results			& $96$ && $71,488$ 	&& $\mathbf{62,475}$	&& $65,442$ 			&& $71,319$  \\
			Average ($n \leq 4$) 	& $51$ && $9,675$ 	&& $7,073$ 				&& $41,300$ 			&& $\mathbf{5,772}$ \\
			Average ($n > 4$)		& $45$ && $141,753$ && $125,417$ 			&& $\mathbf{93,714}$ 	&& $145,733$  \\
			Average (convex) 		& $30$ && $55,320$ 	&& $45,520$ 			&& $42,326$ 			&& $\mathbf{8,950}$  \\
			Average (non-convex) 	& $66$ && $78,837$ 	&& $\mathbf{70,182}$  	&& $75,949$ 			&& $99,669$  \\
			Average (uni-modal)		& $15$ && $28,478$ 	&& $\mathbf{12,624}$ 	&& $23,300$ 			&& $25,796$ \\
			Average (multi-modal)	& $81$ && $81,183$ 	&& $\mathbf{73,979}$  	&& $75,398$ 			&& $81,825$  \\
			Median results			& $96$ && $1,848$ 	&& $960$ 				&& $2,042$ 				&& $\mathbf{775}$  \\
			\bottomrule
			
	\end{tabular}}
	\label{tab:aggresive}
\end{table}

\subsubsection{Investigating the impact of the solution on the boundary}
\label{ssec:boundaries}

In \Cref{sssec:diagonal-trisection}, we stressed that the diagonal trisection strategy (1-DTDV) is especially appealing on problems where the solution lies on the boundary of the feasible region.
We have carried out an additional experimental study presented here to investigate this.
Note that the selection strategy was fixed (IO), and only the influence on the performance of four partitioning strategies was investigated.

Out of the 46 box-constrained unique (excluding dimensionality variations) test problems from \directlib, only the Deb02 problem has a solution on the boundary.
More precisely, the solution lies on the feasible region's vertex, making this situation particularly favorable to the 1-DTDV strategy.
Experimental results using four \direct-type variations on the \textit{Deb02} test problem (with varying dimensionality $n$) are given in the upper part of \Cref{tab:boundaries}.
Since the solution is at the vertex, independently of the dimension, the 1-DTDV-IO algorithm found the solution in the initialization step and took only two objective function evaluations.
The other three \direct-type algorithms required significantly more function evaluations until the stopping condition was satisfied.

In order to carry out a more detailed investigation, test problems with solution coordinates lying on the boundary were artificially constructed.
For this purpose, ten variations of $10$-dimensional \textit{Levy} and five variations of $5$-dimensional \textit{D. Price} test problems with perturbed feasible regions were created.
Let us first consider the case of the \textit{Levy} function.
The original feasible region is $D=[-5, 5]^{10}$ and the solution point $\mathbf{x}^* = (1,1,\dots,1)$.
Thus none of the solution coordinates are on the boundary of the permissible area.
On the original \textit{Levy} problem, the 1-DTDV-IO algorithm performed significantly worse than any of the other three tested \direct-type counterparts (see the first row in the middle part of \Cref{tab:boundaries}).

However, the situation changes completely when \textit{Levy} variations with the perturbed feasible region are considered.
We artiﬁcially reconstructed the original domain so that an increasing number of solution coordinates are located on the boundary---the column ``Feasible region'' in \Cref{tab:boundaries} specifies the modified feasible region.
For example, for the first \textit{Levy} variation (\textit{Levy}1), the modified feasible region coincides with the original one (i.e., $D^m_1 = D$) apart from the coordinate $x_1$, whose new domain is $x_1 \in [-5, 1]$ (original was $x_1 \in [-5, 5]$).
Therefore, for the \textit{Levy}1 problem, one of its solution coordinates $(x_1)$ is on the boundary of the modified feasible region.
The other nine \textit{Levy} function variations are again obtained by substituting only one but the following coordinate compared to the previous \textit{Levy} variation, as shown in \Cref{tab:boundaries}.
The $\nu$ value (the third column in \Cref{tab:boundaries}) indicates the number of solution coordinates projected onto the boundary.
From the obtained results presented in \Cref{tab:boundaries}, observe that the more coordinates of the solution are located on the boundary, the fewer objective function evaluations the 1-DTDV-IO algorithm needs to ﬁnd the solution.
When $\nu \geq 8$ (i.e., at least eight out of 10 coordinates lie on the boundary), the 1-DTDV-IO algorithm outperformed other approaches.
However, the situation is the opposite of the other three partitioning schemes based on \direct-type algorithms.
For almost all cases, they required more function evaluations when $\nu$ increases.

In the last investigation, five variations of \textit{D. Price} test function (original $D = [-10, 10]^5$ and the solution point lying close to the center-point of the domain) were considered.
For \textit{D. Price}, we perturbed the feasible region by the same strategy as for the \textit{Levy} test problem.
Same as before, when most of the solution coordinates are located on the boundaries of the domain ($\nu \geq 4$), the 1-DTDV-IO algorithm is the most efficient.
However, unlike for the \textit{Levy} function, the presence of the solution coordinates on the boundary did not worsen the other \direct-type algorithms but rather improved.
These results show that the performance of center-based partitioning techniques will not necessarily worsen when at least part of the solution coordinates are on the boundary.

\begin{table}[ht]
	\caption{The number of function evaluations required for four different \direct-type algorithms to find optimal solution lying on the boundary}
	\resizebox{\textwidth}{!}{
		\begin{tabular}[tb]{@{\extracolsep{\fill}}lcclrrrr}
			\toprule
			Label & $n$ & $\nu$ & \multicolumn{1}{l}{Feasible region} & N-DTC-IO & 1-DTC-IO & 1-DBDP-IO & 1-DTDV-IO  \\
			\midrule
			\textit{Deb02} & $2$ & $2$ & $D =[0, 1]^n$  & $77$ & $75$ & $100$ & $2$ \\
			\textit{Deb02} & $4$ & $4$ & $D =[0, 1]^n$  & $199$ & $145$ & $220$ & $2$ \\
			\textit{Deb02} & $8$ & $8$ & $D =[0, 1]^n$  & $653$ & $329$ & $494$ & $2$ \\
			\textit{Deb02} & $16$ & $16$ & $D =[0, 1]^n$  & $3,827$ & $1,279$ & $2,368$ & $2$ \\
			\midrule
			\textit{Levy} & $10$ & $0$ & $D = [-5, 5]^n$  & $2,589$ & $919$ & $1,496$ & $141,999$ \\
			\textit{Levy}$1$ & $10$ & $1$ & $D^{\rm m}_1 = D$ and $x_1 \in [-5, 1]$  & $2,847$ & $973$ & $1,326$ & $24,439$ \\
			\textit{Levy}$2$ & $10$ & $2$ & $D^{\rm m}_2 = D^{\rm m}_1$ and $x_2 \in [1, 5]$   & $3,221$ & $1,033$ & $1,386$ & $16,291$ \\
			\textit{Levy}$3$ & $10$ & $3$ & $D^{\rm m}_3 = D^{\rm m}_2$ and $x_3 \in [-10, 1]$ & $3,447$ & $1,079$ & $2,002$ & $13,625$ \\
			\textit{Levy}$4$ & $10$ & $4$ & $D^{\rm m}_4 = D^{\rm m}_3$ and $x_4 \in [1, 10]$  & $3,919$ & $1,119$ & $2,048$ & $14,166$ \\
			\textit{Levy}$5$ & $10$ & $5$ & $D^{\rm m}_5 = D^{\rm m}_4$ and $x_5 \in [-2, 1]$  & $4,091$ & $1,195$ & $2,116$ & $10,927$ \\
			\textit{Levy}$6$ & $10$ & $6$ & $D^{\rm m}_6 = D^{\rm m}_5$ and $x_6 \in [1, 4]$   & $4,483$ & $1,287$ & $2,316$ & $5,416$ \\
			\textit{Levy}$7$ & $10$ & $7$ & $D^{\rm m}_7 = D^{\rm m}_6$ and $x_7 \in [-7, 1]$  & $5,215$ & $2,193$ & $2,484$ & $3,069$ \\
			\textit{Levy}$8$ & $10$ & $8$ & $D^{\rm m}_8 = D^{\rm m}_7$ and $x_8 \in [1, 15]$  & $5,487$ & $2,579$ & $3,174$ & $2,494$ \\
			\textit{Levy}$9$ & $10$ & $9$ & $D^{\rm m}_9 = D^{\rm m}_8$ and $x_9 \in [-13, 1]$ & $6,299$ & $6,581$ & $3,518$ & $547$ \\
			\textit{Levy}$10$ & $10$ & $10$ & $D^{\rm m}_{10} = D^{\rm m}_{9}$ and $x_{10} \in [1,10]$ & $6,487$ & $2,487$ & $3,572$ & $551$ \\
			\midrule
			\textit{D. Price} & $5$ & $0$ & $D = [-10, 10]^n$  & $22,465$ & $20,791$ & $4,060$ & $134,011$ \\
			\textit{D. Price}$1$ & $5$ & $1$ & $D^{\rm m}_1 = D$ and $x_1 \in [-19, 1]$   & $18,245$ & $18,707$ & $2,930$ & $16,089$ \\
			\textit{D. Price}$2$ & $5$ & $2$ & $D^{\rm m}_2 = D^{\rm m}_1$ and $x_2 \in [0.7071, 21]$   & $4,975$ & $1,455$ & $1,322$ & $3,434$ \\
			\textit{D. Price}$3$ & $5$ & $3$ & $D^{\rm m}_3 = D^{\rm m}_2$ and $x_3 \in [-19, 0.5946]$ & $7,709$ & $3,845$ & $1,610$ & $2,759$ \\
			\textit{D. Price}$4$ & $5$ & $4$ & $D^{\rm m}_4 = D^{\rm m}_3$ and $x_4 \in [0.5452, 21]$  & $3,989$ & $1,315$ & $2,280$ & $728$ \\
			\textit{D. Price}$5$ & $5$ & $5$ & $D^{\rm m}_5 = D^{\rm m}_4$ and $x_5 \in [-19, 0.5221]$  & $3,247$ & $1,443$ & $2,396$ & $565$ \\
			
			\bottomrule
			\multicolumn{8}{l}{$\nu$ -- the number of solution coordinates lying on the boundary} \\
	\end{tabular}}
	\label{tab:boundaries}
\end{table}

\subsubsection{Investigating the impact of the domain perturbation}
\label{ssec:perturbation}

By investigating the impact of different partitioning strategies, we observed several situations where one method (or partitioning scheme) dominates the others significantly.
However, sometimes one method may be lucky because the partitioning approach in the initial steps naturally samples near the solution.
In such situations, the location of the solution may favor one partitioning scheme over another.
In this section, we explore whether such dominance is robust to slight perturbations of the domain.

Initially, we identified test problems for which a particular partitioning scheme (regardless of the selection strategy) had a clear dominance, possibly due to the conveniently defined variable bounds.
Out of the $46$ \directlib{} unique box-constrained test problems, the dominance of a particular partitioning scheme was identified for eight of them.
We made domain perturbations for all eight problems keeping the same optimal solution.
First, if any partitioning scheme-based \direct-type algorithm can find the solution in the initialization step, the original domain $(D)$ has been shifted by $22.5 \%$ to the right side.
In some cases, we have made further perturbations that none of the partitioning schemes would sample on initial iterations close to the solution (the column ``Feasible region''  in \Cref{tab:perturb} specifies the original and perturbed domains).


The obtained experimental results revealing the impact of the domain perturbation of twelve \direct-type variants on eight selected original and perturbed test problems with varying dimensionality are given in \Cref{tab:perturb}.
None of the partitioning schemes proved to be robust.
Perturbation of the bounds may significantly reduce the dominance of a particular partitioning scheme.
Quite often, the previously dominant approach may not solve the perturbed problem within the given budget of function evaluations at all.
For example, the 1-DTC partitioning scheme-based \direct-type algorithms undoubtedly dominate original \textit{Rastrigin} and \textit{Griewank} test problems.
However, for perturbed variations of \textit{Rastrigin} and \textit{Griewank}, the same 1-DTC partitioning scheme-based \direct-type algorithms could not be able to find a solution for most of these cases.

Another obvious example is the \textit{Rosenbrock} test problem case.
The 1-DTDV partitioning scheme-based \direct-type algorithms are most efficient when the problem is considered in the original domain $D$.
However, after the $D$'s perturbation, the 1-DTDV partitioning scheme proved to be very inefficient.
Also, in some cases, the bounds' perturbation helped other algorithms perform significantly better than on the original domain.
Examples of such problems are \textit{Styblinski-Tang}, \textit{Easom}, and \textit{Power Sum}.
Furthermore, only for the \textit{Schwefel} test problem, the same two partitioning schemes (1-DBDP and 1-DTDV) based \direct-type algorithms remained the most efficient in the original and perturbed domains.

\begin{sidewaystable}
	\caption{Experimental results of twelve \direct-type variants on 8 selected original and perturbed test problems with varying dimensionality}
	\resizebox{\textwidth}{!}{
		\begin{tabular}[tb]{@{\extracolsep{\fill}}lccrrrrrrrrrrrr}
			\toprule
			\multirow{1}{*}{Label} & \multirow{1}{*}{$n$} & \multirow{1}{*}{Feasible region} & N-DTC-IA & 1-DTC-IA & 1-DBDP-IA & 1-DTDV-IA & N-DTC-IO & 1-DTC-IO & 1-DBDP-IO & 1-DTDV-IO & N-DTC-GL & 1-DTC-GL & 1-DBDP-GL & 1-DTDV-GL  \\
			\midrule
			\textit{Alpine} & $5$ & $[0, 10]^n$ & $61,485$ & $10,343$ & $\mathbf{714}$ & {\color{red} $>10^6$}  & $2,033$ & $1,231$ & $\mathbf{168}$ & $460,445$  & $1,287$ & $685$ & $\mathbf{320}$ & $27,074$	\\
			& $10$ & $[0, 10]^n$ & {\color{red} $>10^6$} & {\color{red} $>10^6$} & $\mathbf{173,912}$ & {\color{red} $>10^6$} & {\color{red} $>10^6$} & {\color{red} $>10^6$} & {\color{red} $>10^6$} & {\color{red} $>10^6$} & $61,209$ & $\mathbf{4,063}$ & $7,646$ & $173,948$  \\
			\textit{Perturbed Alpine} & $5$ & $[\sqrt[i]{2}, 8+\sqrt[i]{2}]^n$ & $61,485$   & $10,343$  & $\mathbf{2,678}$ & $10,941$    & $2,209$ & $\mathbf{1,403}$ & $5,790$ & $3,920$ & $1,479$ & $\mathbf{853}$ & $3,462$ & $19,417$  \\
			& $10$ & $[\sqrt[i]{2}, 8+\sqrt[i]{2}]^n$ & {\color{red} $>10^6$} & {\color{red} $>10^6$} & $\mathbf{218,170}$ & {\color{red} $>10^6$} & {\color{red} $>10^6$} & {\color{red} $>10^6$} & {\color{red} $>10^6$} & {\color{red} $>10^6$} & $31,967$ & $\mathbf{7,245}$ & $15,500$ & {\color{red} $>10^6$}  \\
			\midrule
			\textit{Griewank} & $5$  & $[-330, 870]^n$ & $719,985$ & $69,979$ & {\color{red} $>10^6$} & $\mathbf{47,581}$ & {\color{red} $>10^6$} & {\color{red} $>10^6$} & {\color{red} $>10^6$} & {\color{red} $>10^6$} & $524,765$ & $457,207$ & $\mathbf{14,706}$ & {\color{red} $>10^6$}  \\
			& $10$ & $[-330, 870]^n$ & {\color{red} $>10^6$} & $\mathbf{8,799}$ & $20,534$ & {\color{red} $>10^6$} & {\color{red} $>10^6$} & {\color{red} $>10^6$} & {\color{red} $>10^6$} & {\color{red} $>10^6$} & {\color{red} $>10^6$} & $\mathbf{14,593}$ & $204,940$ & {\color{red} $>10^6$}  \\
			\textit{Perturbed Griewank} & $5$  & $\left[-\sqrt{600i}, \dfrac{600}{\sqrt{i}}\right]^n$ & {\color{red} $>10^6$} & {\color{red} $>10^6$} & {\color{red} $>10^6$} & {\color{red} $>10^6$} & {\color{red} $>10^6$} & {\color{red} $>10^6$} & {\color{red} $>10^6$} & {\color{red} $>10^6$} & $187,811$ & $\mathbf{102,961}$ & $175,806$ & $353,028$   \\
			& $10$ & $\left[-\sqrt{600i}, \dfrac{600}{\sqrt{i}} \right]^n$ & {\color{red} $>10^6$} & {\color{red} $>10^6$} & {\color{red} $>10^6$} & {\color{red} $>10^6$} & {\color{red} $>10^6$} & {\color{red} $>10^6$} & {\color{red} $>10^6$} & {\color{red} $>10^6$} & {\color{red} $>10^6$} & {\color{red} $>10^6$} & {\color{red} $>10^6$} & {\color{red} $>10^6$}  \\
			\midrule
			\textit{Styblinski-Tang} & $5$  & $[-5, 5]^n$ & $4,941$  & $627$   & $\mathbf{164}$ & $243,882$ & $539$   & {\color{red} $>10^6$} & $\mathbf{78}$  & $7,077$ & $1,779$  & $\mathbf{865}$   & $192$ & $682$     \\
			& $10$ & $[-5, 5]^n$ & $68,025$ & $2,631$ & $\mathbf{714}$ & {\color{red} $>10^6$}   & $9,785$ & {\color{red} $>10^6$} & $\mathbf{180}$ & {\color{red} $>10^6$} & $11,347$ & $3,237$ & $\mathbf{784}$ & $5,248$  \\
			\textit{Perturbed Styblinski-Tang} & $5$  & $[-5, 5 + \sqrt[i]{3}]^n$ & $3,919$ & $\mathbf{533}$ & $1,056$ & $66,810$ & $395$ & $\mathbf{273}$ & $278$ & $4,098$ & $1,659$ & $\mathbf{841}$ & $1,654$ & $36,655$  \\
			& $10$ & $[-5, 5 + \sqrt[i]{3}]^n$ & $68,025$ & $\mathbf{2,151}$ & $6,736$ & {\color{red} $>10^6$} & $2,917$ & $\mathbf{829}$ & $1,368$ & {\color{red} $>10^6$} & $13,431$ & $\mathbf{3,447}$ & $11,386$ & {\color{red} $>10^6$}  \\
			\midrule
			\textit{Easom} & $2$ & $[-100, 100]^n$ & $433,031$   & $429,743$ & $\mathbf{444}$ & $20,316$ & $7,581$ & $6,619$ & $\mathbf{322}$ & $6,651$ & $451$ & $\mathbf{321}$ & $544$ & $348$  \\
			\textit{Perturbed Easom} & $2$ & $\left[ \dfrac{-100}{i+1}, 100i \right]^n$ & $214,331$ & $429,743$ & $\mathbf{71,392}$ & $177,258$ & $\mathbf{3,689}$ & $6,659$ & $18,462$ & $10,078$ & $475$ & $393$ & $924$ & $\mathbf{376}$  \\
			\midrule
			\textit{Power Sum} & $4$  & $[0.9, 4.9]^n$ & {\color{red} $>10^6$} & {\color{red} $>10^6$} & $\mathbf{13,828}$ & {\color{red} $>10^6$} & $321,595$ & $144,385$ & $\mathbf{4,790}$ & {\color{red} $>10^6$} & $69,327$ & $77,353$ & $\mathbf{14,214}$ & $37,012$  \\
			\textit{Perturbed Power Sum} & $4$  & $[1, 5 + \sqrt[i]{2}]^n$ & $502,981$ & $176,843$ & $78,746$ & $\mathbf{59,390}$ & $67,959$ & $25,453$ & $17,930$ & $\mathbf{12,219}$ & $152,083$ & $70,745$ & $\mathbf{12,494}$ & $40,753$  \\
			\midrule
			\textit{Rastrigin} & $5$  & $[-2.75, 7.25]^n$ & $8,703$   & $\mathbf{1,487}$ & $314,712$ 			   & {\color{red} $>10^6$} & $597$   & $\mathbf{453}$   & $38,714$              & $112,597$             & $2,721$ & $\mathbf{1,895}$ & $19,642$ & $7,345$  \\
			& $10$ & $[-2.75, 7.25]^n$ & $143,755$ & $\mathbf{7,215}$ & {\color{red} $>10^6$} & {\color{red} $>10^6$} & $4,299$ & $\mathbf{1,551}$ & {\color{red} $>10^6$} & {\color{red} $>10^6$} & $22,971$ & $\mathbf{8,105}$ & {\color{red} $>10^6$} & $140,756$  \\
			\textit{Perturbed Rastrigin} & $5$  & $[-5\sqrt[i]{2}, 7+\sqrt[i]{2}]^n$ & {\color{red} $>10^6$} & {\color{red} $>10^6$} & {\color{red} $>10^6$} & {\color{red} $>10^6$} & {\color{red} $>10^6$} & $\mathbf{567,269}$ & $694,812$ & {\color{red} $>10^6$} & $73,727$ & $24,119$ & $\mathbf{16,440}$ & $90,134$  \\
			& $10$ & $[-5\sqrt[i]{2}, 7+\sqrt[i]{2}]^n$ & {\color{red} $>10^6$} & {\color{red} $>10^6$} & {\color{red} $>10^6$} & {\color{red} $>10^6$} & {\color{red} $>10^6$} & {\color{red} $>10^6$} & {\color{red} $>10^6$} & {\color{red} $>10^6$} & $\mathbf{661,971}$ & {\color{red} $>10^6$} & {\color{red} $>10^6$} & {\color{red} $>10^6$}  \\
			\midrule
			\textit{Rosenbrock} & $5$  & $[-5, 10]^n$  & $73,485$  & $26,325$ & $3,208$ & $\mathbf{1,471}$  & $15,577$ & $1,889$ & $1,494$ & $\mathbf{916}$ & $26,891$  & $15,695$  & $5,110$  & $\mathbf{1,568}$ \\
			& $10$ & $[-5, 10]^n$ & $297,755$ & {\color{red} $>10^6$}  & $13,366$ & $\mathbf{4,541}$ & $71,021$ & {\color{red} $>10^6$} & $4,590$ & $\mathbf{2,091}$  & $104,643$ & $171,019$ & $22,194$ & $\mathbf{5,759}$ \\
			\textit{Perturbed Rosenbrock} & $5$  & $\left[-\dfrac{5}{\sqrt{i}}, 10\sqrt{i}\right]^n$ & $434,985$ & $385,979$ & $\mathbf{291,612}$ & {\color{red} $>10^6$} & $55,693$ & $\mathbf{21,363}$ & $101,508$ & {\color{red} $>10^6$} & $27,763$ & $\mathbf{7,795}$ & $33,056$ & $16,971$  \\
			& $10$ & $\left[-\dfrac{5}{\sqrt{i}}, 10\sqrt{i}\right]^n$ & {\color{red} $>10^6$} & {\color{red} $>10^6$} & {\color{red} $>10^6$} & {\color{red} $>10^6$} & {\color{red} $>10^6$} & {\color{red} $>10^6$} & {\color{red} $>10^6$} & {\color{red} $>10^6$} & $383,081$ & $\mathbf{185,325}$ & $316,392$ & $432,903$  \\
			\midrule
			\textit{Schwefel}  & $5$  & $[-500, 500]^n$ & {\color{red} $>10^6$} & $368,479	$ & $\mathbf{7,566}$ & {\color{red} $>10^6$} & $74,989	$ & $16,767$ & $\mathbf{1,070}$ & $9,561$ & $768,549	$ & $49,247$ & $\mathbf{4,842}$ & $109,746$  \\
			& $10$ & $[-500, 500]^n$ & {\color{red} $>10^6$} & {\color{red} $>10^6$} & $\mathbf{817,512}$ & {\color{red} $>10^6$} & {\color{red} $>10^6$} & {\color{red} $>10^6$} & $\mathbf{57,736}$ & {\color{red} $>10^6$} & {\color{red} $>10^6$} & {\color{red} $>10^6$} & $\mathbf{33,522}$ & {\color{red} $>10^6$}  \\
			\textit{Perturbed Schwefel} & $5$  & $\left[-500 + \dfrac{100}{\sqrt{i}}, 500  - \dfrac{40}{\sqrt{i}} \right]^n$ & {\color{red} $>10^6$} & $458,979$ & $19,972$ & $\mathbf{2,135}$ & $80,295$ & $35,091$ & $84,096$ & $\mathbf{33,622}$ & $336,581$ & $65,329$ & $9,548$ & $\mathbf{1,580}$  \\
			& $10$ & $\left[-500 + \dfrac{100}{\sqrt{i}}, 500  - \dfrac{40}{\sqrt{i}} \right]^n$ & {\color{red} $>10^6$} & {\color{red} $>10^6$} & {\color{red} $>10^6$} & {\color{red} $>10^6$} & {\color{red} $>10^6$} & {\color{red} $>10^6$} & {\color{red} $>10^6$} & {\color{red} $>10^6$} & {\color{red} $>10^6$} & {\color{red} $>10^6$} & $99,824$ & $336,983$  \\
			\bottomrule
			\multicolumn{7}{l}{$i = 1,...,n$ -- indexes used for variable bounds}
	\end{tabular}}
	\label{tab:perturb}
\end{sidewaystable}

\subsection{Investigation of different selection schemes}
\label{ssec:selection-impact}

Here, the efficiency of new \direct-type algorithms is investigated based on the selection scheme.
In \Cref{tab:aggresive}, three row parts corresponds to a different selection approach.
Since each selection strategy was run on the $96$ problems using $4$ different partitioning methods (columns of \Cref{tab:aggresive}), it follows that each selection approach was involved in solving $4 \times 96 = 384$ test problems.
We note that algorithms incorporating a two-step-based Pareto selection scheme (GL) combined with any partitioning strategy, on average, deliver the best results.
All algorithmic variants based on the GL selection scheme did not solve ($19/384$).
In contrast, the IO and IA selection schemes based algorithms failed to solve ($63/384$) and ($55/384$) cases accordingly.
This leads to a much better average performance of \direct-type algorithms based on the GL selection scheme.
In total, GL selection scheme-based algorithms required approximately $60\%$ and $62\%$ fever function evaluations compared with the IO and IA counterparts.
The GL selection scheme's most significant advantage can be seen in solving higher-dimensional $(n > 4)$ test problems.
In total, GL selection scheme-based algorithms solving ($n > 4$) test instances required approximately $72\%$ fever function evaluations compared with the IO or IA selection scheme-based counterparts.

The IO selection scheme seems the most suitable for simpler optimization problems (low dimensional, uni-modal, and problems with a few minima).
Of 63 failed problems using the IO selection scheme, 43 were extremely hard, i.e., multi-modal, sharply peaked, and multi-variable (e.g., $n \ge 10$).
The GL selection strategy-based variants usually select more regions to subdivide.
Therefore, they suffer for these more straightforward optimization problems.
However, the GL scheme ensured, that on average, all \direct-type variants converged in significantly fewer function evaluations on complex multi-modal test problems.

\begin{figure}
	\resizebox{\textwidth}{!}{
		\begin{tikzpicture}
			\begin{axis}[
				legend pos=north west,
				title  = {Operational characteristics},
				xlabel = {Function evaluations},
				xmode=log,
				ymin=-0.02,ymax=1.02,
				ytick distance=0.1,
				xmode=log,
				xmin=10,
				xmax=1000000,
				xtick distance=10,
				ylabel = {Proportion of solved problems},
				ylabel style={yshift=-0.5em},
				legend style={font=\tiny,xshift=-0.5em},
				legend cell align={left},
				legend columns=1,
				height=0.75\textwidth,width=\textwidth,
				every axis plot/.append style={very thick},
				]
				\addplot[mark=*,black,mark options={scale=1.5, fill=princetonorange}, only marks,line width=0.75pt] coordinates {(0.1,0.1)} ;
				\label{p1}
				\addplot[mark=square*,black,mark options={scale=1.5, fill=yellow}, only marks,line width=0.75pt] coordinates {(0.1,0.1)} ;
				\label{p2}
				\addplot[mark=diamond*,black,mark options={scale=1.5, fill=sienna}, only marks,line width=0.75pt] coordinates {(0.1,0.1)} ;
				\label{p3}
				\addplot[mark=triangle*,black,mark options={scale=1.5, fill=psychedelicpurple}, only marks,line width=0.75pt] coordinates {(0.1,0.1)} ;
				\label{p4}
				\node [draw,fill=white] at (rel axis cs: 0.85,0.14) {\shortstack[l]{
						{\scriptsize \textbf{Partitioning  schemes}} \\
						\ref{p1} {\scriptsize N-DTC} \\
						\ref{p2} {\scriptsize 1-DTC} \\
						\ref{p3} {\scriptsize 1-DBDP} \\
						\ref{p4} {\scriptsize 1-DTDV}}};
				
				\addplot[postaction={decoration={markings,mark=between positions 0 and 1 step 0.06 with {\node[circle,draw=black,fill=princetonorange,inner sep=1.5pt,solid] {};}},decorate,},sandstorm,line width=0.75pt,densely dashed] table[x=T,y=DDA] {data/Overallass.txt};
				\addplot[postaction={decoration={markings,mark=between positions 0 and 1 step 0.07 with {\node[mark=square,draw=black,fill=yellow,inner sep=1.5pt,solid] {};}},decorate,},sandstorm,line width=0.75pt,densely dashed] table[x=T,y=DRA] {data/Overallass.txt};
				\addplot[postaction={decoration={markings,mark=between positions 0 and 1 step 0.08 with 	{\node[diamond,draw=black,fill=sienna,inner sep=1.5pt,solid] {};}},decorate,},sandstorm,line width=0.75pt,densely dashed] table[x=T,y=BIA] {data/Overallass.txt};
				\addplot[postaction={decoration={markings,mark=between positions 0 and 1 step 0.09 with {\node[regular 	polygon,regular polygon sides=3,draw=black,fill=psychedelicpurple,inner sep=1pt,solid] {};}},decorate,},sandstorm,line width=0.75pt,densely dashed] table[x=T,y=ADA] {data/Overallass.txt};
				
				\addplot[postaction={decoration={markings,mark=between positions 0 and 1 step 0.1 with {\node[circle,draw=black,fill=princetonorange,inner sep=1.5pt] {};}},decorate,},blue,line width=0.75pt] table[x=T,y=DDO] {data/Overallass.txt};
				\addplot[postaction={decoration={markings,mark=between positions 0 and 1 step 0.11 with {\node[mark=square,draw=black,fill=yellow,inner sep=1.5pt] {};}},decorate,},blue,line width=0.75pt] table[x=T,y=DRO] {data/Overallass.txt};
				\addplot[postaction={decoration={markings,mark=between positions 0 and 1 step 0.12 with 	{\node[diamond,draw=black,fill=sienna,inner sep=1.5pt] {};}},decorate,},blue,line width=0.75pt] table[x=T,y=BIO] {data/Overallass.txt};
				\addplot[postaction={decoration={markings,mark=between positions 0 and 1 step 0.13 with {\node[regular 	polygon,regular polygon sides=3,draw=black,fill=psychedelicpurple,inner sep=1pt] {};}},decorate,},blue,line width=0.75pt] table[x=T,y=ADO] {data/Overallass.txt};
				
				\addplot[postaction={decoration={markings,mark=between positions 0 and 1 step 0.14 with {\node[circle,draw=black,fill=princetonorange,inner sep=1.5pt,solid] {};}},decorate,},onyx,line width=0.75pt,densely dotted] table[x=T,y=DDG] {data/Overallass.txt};
				\addplot[postaction={decoration={markings,mark=between positions 0 and 1 step 0.15 with {\node[mark=square,draw=black,fill=yellow,inner sep=1.5pt,solid] {};}},decorate,},onyx,line width=0.75pt,densely dotted] table[x=T,y=DRG] {data/Overallass.txt};
				\addplot[postaction={decoration={markings,mark=between positions 0 and 1 step 0.16 with 	{\node[diamond,draw=black,fill=sienna,inner sep=1.5pt,solid] {};}},decorate,},onyx,line width=0.75pt,densely dotted] table[x=T,y=BIG] {data/Overallass.txt};
				\addplot[postaction={decoration={markings,mark=between positions 0 and 1 step 0.17 with {\node[regular 	polygon,regular polygon sides=3,draw=black,fill=psychedelicpurple,inner sep=1pt,solid] {};}},decorate,},onyx,line width=0.75pt,densely dotted] table[x=T,y=ADG] {data/Overallass.txt};
				
				\addplot[blue,line width=0.75pt] coordinates {(0.1,0.1)} ;
				\label{p5}
				\addplot[sandstorm,line width=0.75pt,densely dashed] coordinates {(0.1,0.1)} ;
				\label{p6}
				\addplot[onyx,line width=0.75pt,densely dotted] coordinates {(0.1,0.1)} ;
				\label{p7}
				\node [draw,fill=white] at (rel axis cs: 0.2,0.88) {\shortstack[l]{
						{\scriptsize \textbf{POH selection schemes}} \\
						\ref{p5} {\scriptsize Improved Original (IO)} \\
						\ref{p6} {\scriptsize Improved Aggressive (IA)} \\
						\ref{p7} {\scriptsize Global-Local (GL) Pareto}}};
			\end{axis}
	\end{tikzpicture}}
	
	\caption{Operational characteristics for all twelve \direct-type algorithmic variations on \directlib{} test problems}
	\label{fig:perf-l1}
\end{figure}
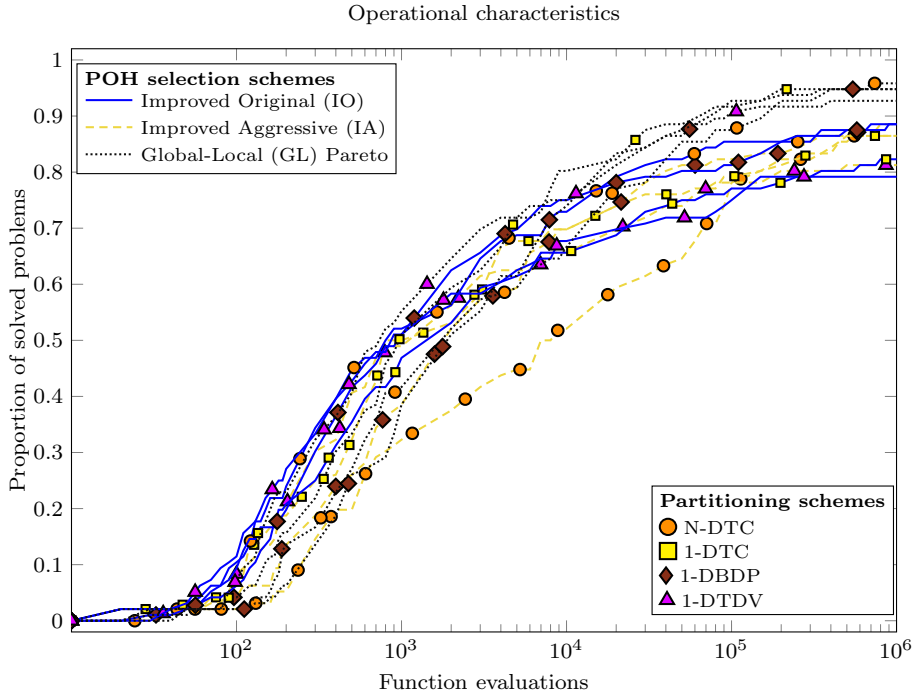

\begin{figure}
	\resizebox{\textwidth}{!}{
		\begin{tikzpicture}
			\begin{axis}[
				legend pos=north west,
				title  = {Operational characteristics},
				xlabel = {Function evaluations},
				xmode=log,
				ymin=-0.02,ymax=1.02,
				ytick distance=0.1,
				xmode=log,
				xmin=10,
				xmax=1000000,
				xtick distance=10,
				ylabel = {Proportion of solved problems},
				ylabel style={yshift=-0.5em},
				legend style={font=\tiny,xshift=-0.5em},
				legend cell align={left},
				legend columns=1,
				height=0.75\textwidth,width=\textwidth,
				every axis plot/.append style={very thick},
				]
				\node [draw,fill=white] at (rel axis cs: 0.85,0.14) {\shortstack[l]{
						{\scriptsize \textbf{Partitioning  schemes}} \\
						\ref{p1} {\scriptsize N-DTC} \\
						\ref{p2} {\scriptsize 1-DTC} \\
						\ref{p3} {\scriptsize 1-DBDP} \\
						\ref{p4} {\scriptsize 1-DTDV}}};
				
				\addplot[postaction={decoration={markings,mark=between positions 0 and 1 step 0.06 with {\node[circle,draw=black,fill=princetonorange,inner sep=1.5pt,solid] {};}},decorate,},sandstorm,line width=0.75pt,densely dashed] table[x=T,y=DDA] {data/Overallassa.txt};
				\addplot[postaction={decoration={markings,mark=between positions 0 and 1 step 0.07 with {\node[mark=square,draw=black,fill=yellow,inner sep=1.5pt,solid] {};}},decorate,},sandstorm,line width=0.75pt,densely dashed] table[x=T,y=DRA] {data/Overallassa.txt};
				\addplot[postaction={decoration={markings,mark=between positions 0 and 1 step 0.08 with 	{\node[diamond,draw=black,fill=sienna,inner sep=1.5pt,solid] {};}},decorate,},sandstorm,line width=0.75pt,densely dashed] table[x=T,y=BIA] {data/Overallassa.txt};
				\addplot[postaction={decoration={markings,mark=between positions 0 and 1 step 0.09 with {\node[regular 	polygon,regular polygon sides=3,draw=black,fill=psychedelicpurple,inner sep=1pt,solid] {};}},decorate,},sandstorm,line width=0.75pt,densely dashed] table[x=T,y=ADA] {data/Overallassa.txt};
				
				\addplot[postaction={decoration={markings,mark=between positions 0 and 1 step 0.1 with {\node[circle,draw=black,fill=princetonorange,inner sep=1.5pt] {};}},decorate,},blue,line width=0.75pt] table[x=T,y=DDO] {data/Overallassa.txt};
				\addplot[postaction={decoration={markings,mark=between positions 0 and 1 step 0.11 with {\node[mark=square,draw=black,fill=yellow,inner sep=1.5pt] {};}},decorate,},blue,line width=0.75pt] table[x=T,y=DRO] {data/Overallassa.txt};
				\addplot[postaction={decoration={markings,mark=between positions 0 and 1 step 0.12 with 	{\node[diamond,draw=black,fill=sienna,inner sep=1.5pt] {};}},decorate,},blue,line width=0.75pt] table[x=T,y=BIO] {data/Overallassa.txt};
				\addplot[postaction={decoration={markings,mark=between positions 0 and 1 step 0.13 with {\node[regular 	polygon,regular polygon sides=3,draw=black,fill=psychedelicpurple,inner sep=1pt] {};}},decorate,},blue,line width=0.75pt] table[x=T,y=ADO] {data/Overallassa.txt};
				
				\addplot[postaction={decoration={markings,mark=between positions 0 and 1 step 0.14 with {\node[circle,draw=black,fill=princetonorange,inner sep=1.5pt,solid] {};}},decorate,},onyx,line width=0.75pt,densely dotted] table[x=T,y=DDG] {data/Overallassa.txt};
				\addplot[postaction={decoration={markings,mark=between positions 0 and 1 step 0.15 with {\node[mark=square,draw=black,fill=yellow,inner sep=1.5pt,solid] {};}},decorate,},onyx,line width=0.75pt,densely dotted] table[x=T,y=DRG] {data/Overallassa.txt};
				\addplot[postaction={decoration={markings,mark=between positions 0 and 1 step 0.16 with 	{\node[diamond,draw=black,fill=sienna,inner sep=1.5pt,solid] {};}},decorate,},onyx,line width=0.75pt,densely dotted] table[x=T,y=BIG] {data/Overallassa.txt};
				\addplot[postaction={decoration={markings,mark=between positions 0 and 1 step 0.17 with {\node[regular 	polygon,regular polygon sides=3,draw=black,fill=psychedelicpurple,inner sep=1pt,solid] {};}},decorate,},onyx,line width=0.75pt,densely dotted] table[x=T,y=ADG] {data/Overallassa.txt};
				
				\node [draw,fill=white] at (rel axis cs: 0.2,0.88) {\shortstack[l]{
						{\scriptsize \textbf{POH selection schemes}} \\
						\ref{p5} {\scriptsize Improved Original (IO)} \\
						\ref{p6} {\scriptsize Improved Aggressive (IA)} \\
						\ref{p7} {\scriptsize Global-Local (GL) Pareto}}};
			\end{axis}
	\end{tikzpicture}}
	
	\caption{Operational characteristics for all twelve \direct-type algorithmic variations solving higher-dimensional ($n > 4$) multi-modal \directlib{} test problems}
	\label{fig:perf-l1d}
\end{figure}
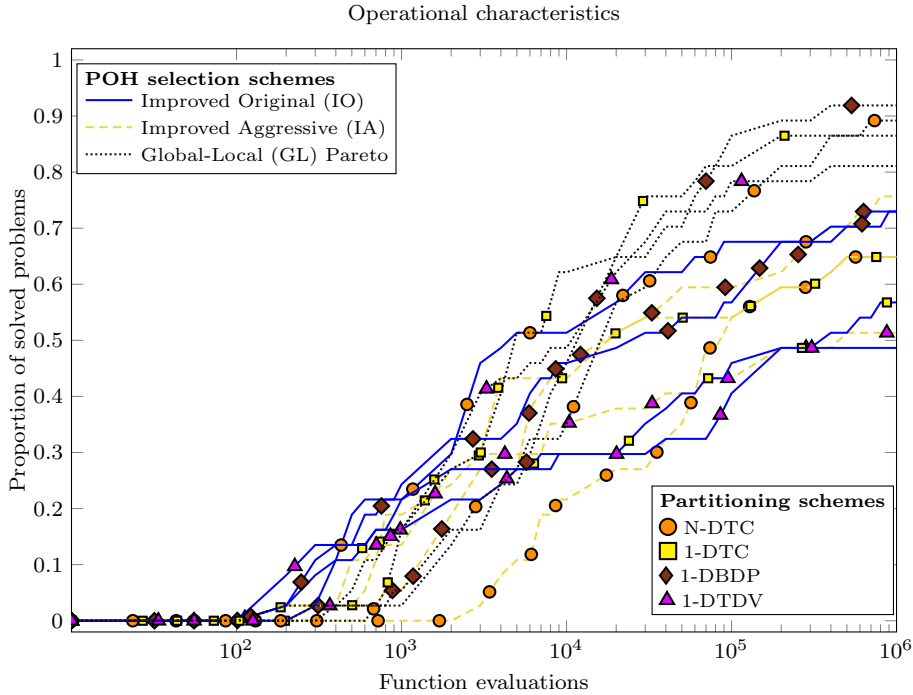

\begin{figure}
	\resizebox{\textwidth}{!}{
		\begin{tikzpicture}
			\begin{axis}[
				legend pos=north west,
				title  = {Operational characteristics},
				xlabel = {Function evaluations},
				xmode=log,
				ymin=-0.02,ymax=1.02,
				ytick distance=0.1,
				xmode=log,
				xmin=10,
				xmax=1000000,
				xtick distance=10,
				ylabel = {Proportion of solved problems},
				ylabel style={yshift=-0.5em},
				legend style={font=\tiny,xshift=-0.5em},
				legend cell align={left},
				legend columns=1,
				height=0.75\textwidth,width=\textwidth,
				every axis plot/.append style={very thick},
				]
				\node [draw,fill=white] at (rel axis cs: 0.85,0.14) {\shortstack[l]{
						{\scriptsize \textbf{Partitioning  schemes}} \\
						\ref{p1} {\scriptsize N-DTC} \\
						\ref{p2} {\scriptsize 1-DTC} \\
						\ref{p3} {\scriptsize 1-DBDP} \\
						\ref{p4} {\scriptsize 1-DTDV}}};
				
				\addplot[postaction={decoration={markings,mark=between positions 0 and 1 step 0.06 with {\node[circle,draw=black,fill=princetonorange,inner sep=1.5pt,solid] {};}},decorate,},sandstorm,line width=0.75pt,densely dashed] table[x=T,y=DDA] {data/Overallassb.txt};
				\addplot[postaction={decoration={markings,mark=between positions 0 and 1 step 0.07 with {\node[mark=square,draw=black,fill=yellow,inner sep=1.5pt,solid] {};}},decorate,},sandstorm,line width=0.75pt,densely dashed] table[x=T,y=DRA] {data/Overallassb.txt};
				\addplot[postaction={decoration={markings,mark=between positions 0 and 1 step 0.08 with 	{\node[diamond,draw=black,fill=sienna,inner sep=1.5pt,solid] {};}},decorate,},sandstorm,line width=0.75pt,densely dashed] table[x=T,y=BIA] {data/Overallassb.txt};
				\addplot[postaction={decoration={markings,mark=between positions 0 and 1 step 0.09 with {\node[regular 	polygon,regular polygon sides=3,draw=black,fill=psychedelicpurple,inner sep=1pt,solid] {};}},decorate,},sandstorm,line width=0.75pt,densely dashed] table[x=T,y=ADA] {data/Overallassb.txt};
				
				\addplot[postaction={decoration={markings,mark=between positions 0 and 1 step 0.1 with {\node[circle,draw=black,fill=princetonorange,inner sep=1.5pt] {};}},decorate,},blue,line width=0.75pt] table[x=T,y=DDO] {data/Overallassb.txt};
				\addplot[postaction={decoration={markings,mark=between positions 0 and 1 step 0.11 with {\node[mark=square,draw=black,fill=yellow,inner sep=1.5pt] {};}},decorate,},blue,line width=0.75pt] table[x=T,y=DRO] {data/Overallassb.txt};
				\addplot[postaction={decoration={markings,mark=between positions 0 and 1 step 0.12 with 	{\node[diamond,draw=black,fill=sienna,inner sep=1.5pt] {};}},decorate,},blue,line width=0.75pt] table[x=T,y=BIO] {data/Overallassb.txt};
				\addplot[postaction={decoration={markings,mark=between positions 0 and 1 step 0.13 with {\node[regular 	polygon,regular polygon sides=3,draw=black,fill=psychedelicpurple,inner sep=1pt] {};}},decorate,},blue,line width=0.75pt] table[x=T,y=ADO] {data/Overallassb.txt};
				
				\addplot[postaction={decoration={markings,mark=between positions 0 and 1 step 0.14 with {\node[circle,draw=black,fill=princetonorange,inner sep=1.5pt,solid] {};}},decorate,},onyx,line width=0.75pt,densely dotted] table[x=T,y=DDG] {data/Overallassb.txt};
				\addplot[postaction={decoration={markings,mark=between positions 0 and 1 step 0.15 with {\node[mark=square,draw=black,fill=yellow,inner sep=1.5pt,solid] {};}},decorate,},onyx,line width=0.75pt,densely dotted] table[x=T,y=DRG] {data/Overallassb.txt};
				\addplot[postaction={decoration={markings,mark=between positions 0 and 1 step 0.16 with 	{\node[diamond,draw=black,fill=sienna,inner sep=1.5pt,solid] {};}},decorate,},onyx,line width=0.75pt,densely dotted] table[x=T,y=BIG] {data/Overallassb.txt};
				\addplot[postaction={decoration={markings,mark=between positions 0 and 1 step 0.17 with {\node[regular 	polygon,regular polygon sides=3,draw=black,fill=psychedelicpurple,inner sep=1pt,solid] {};}},decorate,},onyx,line width=0.75pt,densely dotted] table[x=T,y=ADG] {data/Overallassb.txt};
				
				\node [draw,fill=white] at (rel axis cs: 0.2,0.88) {\shortstack[l]{
						{\scriptsize \textbf{POH selection schemes}} \\
						\ref{p5} {\scriptsize Improved Original (IO)} \\
						\ref{p6} {\scriptsize Improved Aggressive (IA)} \\
						\ref{p7} {\scriptsize Global-Local (GL) Pareto}}};
			\end{axis}
	\end{tikzpicture}}
	
	\caption{Operational characteristics for all twelve \direct-type algorithmic variations solving uni-modal and convex  \directlib{} test problems}
	\label{fig:perf-l1c}
\end{figure}
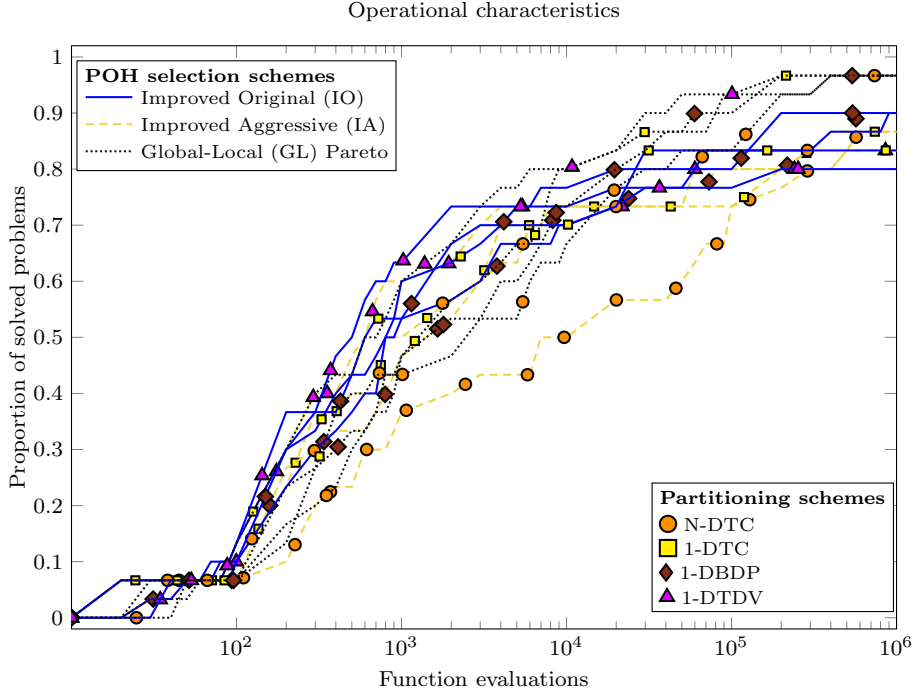

Additionally, the operational characteristics~\cite{Grishagin1978,Strongin2000:book} reported in \Cref{fig:perf-l1,fig:perf-l1d,fig:perf-l1c} show the behavior of all twelve algorithms on different subsets of \directlib{} box-constrained test problems.
Operational characteristics provide the proportion of problems that can be solved within a given budget of function evaluations.
\Cref{fig:perf-l1}, drawn using all $96$ box-constrained problems from \directlib{}, clearly shows that IO selection scheme-based \direct-type algorithms (1-DTDV-IO and N-DTC-IO) dominate for simpler problems.
They solved almost half of the 96 test problems within a small budget of objective function evaluations.
However, as the number of function evaluations increases (as more complex problems are considered), the GL scheme-based algorithms are most efficient.


Operational characteristics in \Cref{fig:perf-l1d} show the behavior of all twelve algorithms on $35$ higher-dimensionality ($n > 4$) multi-modal test problems.
Once again, when a given budget of function evaluations is low (M$_{\rm max} \leq 500$), all IO selection scheme-based variations perform better.
Unfortunately, with such a small function evaluation budget, the algorithms will only solve approximately $15 \%$ of all the test problems.
When the maximal budget of function evaluations increased (M$_{\rm max} \leq 4,000$), only one of the IO selection scheme combinations (N-DTC-IO) maintained the highest efficiency and solved approximately $50 \%$ of all the test cases.
Finally, when the function evaluation budget is higher (M$_{\rm max} \geq 4,000$), GL selection scheme-based variations (1-DTC-GL and 1-DBDP-GL) have the highest efficiency.

Similar tendencies regarding the best-performing selection strategies can be seen in \Cref{fig:perf-l1c}.
Here, the operational characteristics illustrate the behavior of algorithms solving simplest uni-modal and convex optimization test problems.
Among the partitioning strategies, the 1-DTDV scheme looks the most efficient here.

\section{Experimental investigation using GKLS-type test problems}
\label{sec:exp-GKLS}

Additionally, we compare the performance of all twelve \direct-type variants on GKLS-type test problems~\cite{Gaviano2003}.
GKLS-generator allows generating three types (non-differentiable, continuously differentiable, and twice continuously differentiable) of multi-dimensional and multi-extremal optimization test functions with a priori known local and global minima.
The complexity of generated problems is established by setting different values for user-determined parameters: problem dimension $n$, the number of local minima $m$, global minimum value $f^*$, distance $d$ from the global minimizer to the paraboloid vertex, and radius $r$ of the attraction region of the global minimizer.

We use eight different complexity classes (see \Cref{tab:paramet}).
The dimensionality ($n$) and other parameters are set as in~\cite{Paulavicius2014:jogo,Paulavicius2019:eswa}.
Each class consisted of $100$ test instances.
For each dimension $n$, two test classes were considered: the ``simple'' class and the ``hard'' one.
For three- and four-dimensional classes the difficulty is increased by enlarging the distance $d$ from the global minimizer $(\mathbf{x}^*)$ to the paraboloid vertex.
For two and five-dimensional classes this is achieved by decreasing the radius $r$ of the attraction region of the global minimizer.

\begin{table}[ht]
	\centering
	\caption{Description of GKLS-type test classes used in numerical experiments}
	\label{tab:paramet}
	\resizebox{\textwidth}{!}{
		\begin{tabular}{p{1cm}p{2cm}p{1cm}p{1cm}p{1cm}p{1cm}p{1cm}p{1cm}}
			\toprule
			Class & Difficulty & $\Delta$ & $n$ & $f^*$ & $d$ & $r$ & $m$ \\
			\midrule
			$1$ & simple 	& $10^{-4}$ & $2$ & $-1$ & $0.90$ & $0.20$ & $10$  \\
			$2$ & hard 		& $10^{-4}$ & $2$ & $-1$ & $0.90$ & $0.10$ & $10$  \\
			$3$ & simple 	& $10^{-6}$ & $3$ & $-1$ & $0.66$ & $0.20$ & $10$  \\
			$4$ & hard 		& $10^{-6}$ & $3$ & $-1$ & $0.90$ & $0.20$ & $10$  \\
			$5$ & simple 	& $10^{-6}$ & $4$ & $-1$ & $0.66$ & $0.20$ & $10$  \\
			$6$ & hard 		& $10^{-6}$ & $4$ & $-1$ & $0.90$ & $0.20$ & $10$  \\
			$7$ & simple 	& $10^{-7}$ & $5$ & $-1$ & $0.66$ & $0.30$ & $10$  \\
			$8$ & hard 		& $10^{-7}$ & $5$ & $-1$ & $0.66$ & $0.20$ & $10$  \\
			\bottomrule
	\end{tabular}}
\end{table}

The same stopping rule is adopted in these experiments as in~\cite{Paulavicius2014:jogo,Paulavicius2019:eswa}.
The global minimizer $\mathbf{x}^* \in D$ is considered to be found when an algorithm generated a function evaluation point $\mathbf{x}^i \in D^i_k$ such that:
\begin{equation}
	\label{eq:pes}
	\arrowvert x^i_j - x^*_j \arrowvert \leq \sqrt[n]{\Delta}(b_j - a_j), \hspace{1cm} 1 \leq j \leq n,
\end{equation}
where $0 \leq \Delta \leq 1$ is an accuracy coefficient \cite{Sergeyev2006}(see \Cref{tab:paramet} for $\Delta$ parameter values).
In other words, we stop the search when the algorithm has produced a point very close to the known optimum.
In each run, we used the same limit of function evaluations equal to $10^6$.
Note that the stopping rule \eqref{eq:pes} does not require algorithms to find a solution with high accuracy.
Therefore, the Pareto selection enhancing the local search (see \Cref{sssec:two-step-selection}) was disabled.

The experimental results are summarized in \Cref{tab:results4}.
The notation “$>10^6(j)$” indicates that after the maximal number of function evaluations $10^6$, the algorithm under consideration was not able to solve $j$ problems in total.
First, contrary to the previous tendencies, the best results are obtained when \direct-type variants include an improved original selection scheme (IO).
In numerical experiments described in \Cref{sec:experiments}, \direct-type algorithms with integrated IO selection schemes had almost the worst efficiencies.
However, all four \direct-type variants based on the IO selection scheme are promising for GKLS-type test problems.
The least attractive is an improved aggressive selection scheme (IA).
All failed $40$ cases appeared when this selection scheme was combined with three different partitioning strategies (except 1-DBDP).

The best average results (see the upper part of \Cref{tab:results4}) are achieved using 1-DTC-IO (for seven different classes) and 1-DBDP-GL (for one class).
The average number using the 1-DTC-IO algorithm is $8,021$, while the second (1-DBDP-IO) and third-best (1-DTDV-IO) algorithms deliver approximately $19 \%$ ($9,571$) and $47 \%$ ($15,231$) worse overall performances.
Interestingly, the 1-DBDP-GL algorithm, which produced the best overall result in \Cref{sec:experiments}, ranks only fifth as delivered $67 \%$ ($24,746$) worse average results than 1-DTC-IO.

The lowest aggregated median number (the middle part of \Cref{tab:results4}) for all eight classes again is obtained using the same 1-DTC-IO algorithm ($1,427$).
Therefore, the 1-DTC-IO algorithm can solve at least half of GKLS-type problems with the best performance.
In contrast, the second (1-DBDP-GL) and third (1-DBDP-IO) best algorithms delivered approximately $3.71 \%$ ($1,482$) and $3.84 \%$ ($1,484$) worse overall median values.

In the bottom part of \Cref{tab:results4}, the maximal number of function evaluations required to solve test problems within a particular class is reported.
The previously emphasized algorithmic variation 1-DTC-IO is the best for four out of eight classes.
However, solving two of the most complex (``hard'') classes (No = $6$ and $8$) 1-DTDV-IO algorithm seems the most promising.

\begin{sidewaystable}
	\caption{Comparison of twelve \direct-type variants on eight classes of GKLS-type problems}
	\resizebox{\textwidth}{!}{
		\begin{tabular}[tb]{@{\extracolsep{\fill}}c|rrrr|rrrr|rrrr}
			\toprule
			Class & N-DTC-IA & 1-DTC-IA & 1-DBDP-IA & 1-DTDV-IA & N-DTC-IO & 1-DTC-IO & 1-DBDP-IO & 1-DTDV-IO & N-DTC-GL & 1-DTC-GL & 1-DBDP-GL & 1-DTDV-GL \\
			\midrule
			\multicolumn{13}{c}{Average number of function evaluations}\\
			\midrule
			$1$ & $510$      & $289$      & $257$     & $1,148$   & $198$     & $\mathbf{149}$     & $156$    & $360$    & $236$     & $217$     & $185$    & $665$ \\
			$2$ & $5,215$    & $5,491$    & $1,449$   & $6,278$   & $1,068$   & $\mathbf{803}$     & $863$    & $1,458$  & $2,332$   & $3,130$   & $1,245$  & $3,454$ \\
			$3$ & $4,069$    & $1,747$    & $1,252$   & $9,870$   & $1,019$   & $\mathbf{673}$     & $869$    & $1,776$  & $1,264$   & $1,286$   & $825$    & $4,148$ \\
			$4$ & $15,331$   & $9,393$    & $4,032$   & $35,291$  & $2,477$   & $\mathbf{1,543}$   & $1,832$  & $4,539$  & $5,125$   & $4,427$   & $2,640$  & $14,943$ \\
			$5$ & $36,586$   & $20,911$   & $13,179$  & $91,161$  & $7,843$   & $\mathbf{4,591}$   & $8,207$  & $11,806$ & $10,720$  & $10,553$  & $9,199$  & $32,597$ \\
			$6$ & $198,129$  & $157,187$  & $95,541$  & $350,033$ & $26,970$  & $\mathbf{16,899}$  & $23,905$ & $34,020$ & $77,715$  & $68,348$  & $57,867$ & $117,684$ \\
			$7$ & $21,811$   & $17,769$   & $4,956$   & $67,634$  & $6,216$   & $4,419$ & $3,974$  & $11,156$ & $12,262$ & $10,645$  & $\mathbf{3,472}$  & $32,108$ \\
			$8$ & $308,948$  & $226,208$  & $96,621$  & $477,868$ & $67,636$  & $\mathbf{35,091}$  & $36,768$ & $55,625$ & $120,118$ & $99,362$  & $58,278$ & $183,881$ \\
			$1-8$ & $73,825$ & $54,874$   & $27,160$  & $129,910$ & $14,178$  & $\mathbf{8,021}$   & $9,571$  & $15,231$ & $27,783$  & $24,746$  & $16,713$ & $48,685$ \\
			\midrule
			\multicolumn{13}{c}{Median number of function evaluations}\\
			\midrule
			$1$ & $308$      & $207$      & $204$     & $878$     & $119$     & $122$     		 & $\mathbf{111}$    & $328$    & $128$    & $155$      & $130$     & $469$ \\
			$2$ & $4,316$    & $2,375$    & $1,016$   & $5,793$   & $1,058$   & $724$     		 & $\mathbf{673}$    & $1,364$  & $1,979$  & $2,097$    & $855$     & $3,272$ \\
			$3$ & $1,267$    & $949$      & $755$     & $7,813$   & $\mathbf{387}$ & $503$  	 & $488$    	     & $1,597$  & $445$    & $604$      & $461$     & $3,344$ \\
			$4$ & $8,970$    & $3,962$    & $2,366$   & $32,915$  & $1,782$   & $\mathbf{1,075}$ & $1,189$  	     & $4,230$  & $2,463$  & $2,992$    & $1,399$   & $13,753$ \\
			$5$ & $15,523$   & $9,063$    & $5,851$   & $78,151$  & $4,874$   & $\mathbf{2,872}$ & $4,443$  	   	 & $10,761$ & $4,388$  & $7,052$    & $4,202$   & $28,822$ \\
			$6$ & $121,285$  & $76,586$   & $55,130$  & $334,277$ & $15,517$  & $\mathbf{9,237}$ & $15,628$ 	  	 & $32,796$ & $43,458$ & $33,804$   & $29,189$  & $116,530$ \\
			$7$ & $7,514$    & $4,079$    & $2,102$   & $42,804$  & $1,673$   & $2,291$   		 & $2,278$  		 & $10,992$ & $1,533$  & $3,440$    & $\mathbf{1,427}$   & $29,635$ \\
			$8$ & $203,711$  & $128,408$  & $53,291$  & $429,811$ & $43,400$  & $24,327$  		 & $\mathbf{19,967}$ & $47,221$ & $65,892$ & $54,497$   & $27,697$  & $162,411$ \\
			$1-8$ & $6,767$  & $4,188$    & $2,098$   & $26,227$  & $1,644$   & $\mathbf{1,427}$ & $1,487$  		 & $4,599$  & $2,117$  & $3,178$  & $1,482$  & $13,235$  \\
			\midrule
			\multicolumn{13}{c}{Maximal number of function evaluations}\\
			\midrule
			$1$ & $4,777$    & $1,955$    & $1,178$   & $4,811$     & $1,153$   & $\mathbf{655}$     & $840$     & $961$     & $2,031$   & $1,319$   & $1,090$   & $2,502$   \\
			$2$ & $21,841$   & $25,021$   & $11,674$  & $20,089$    & $3,197$   & $\mathbf{2,201}$   & $4,374$   & $3,964$   & $8,431$   & $11,041$  & $8,836$   & $10,269$  \\
			$3$ & $31,291$   & $13,233$   & $8,196$   & $34,607$    & $6,625$   & $\mathbf{3,273}$   & $5,032$   & $4,864$   & $10,723$  & $11,335$  & $5,058$   & $16,789$  \\
			$4$ & $132,121$  & $150,809$  & $22,720$  & $88,327$    & $15,307$  & $9,763$   & $\mathbf{7,806}$   & $10,236$  & $48,371$  & $45,721$  & $14,064$  & $37,514$  \\
			$5$ & $212,339$  & $131,783$  & $116,136$ & $280,015$   & $39,129$  & $\mathbf{18,853}$  & $62,016$  & $35,898$  & $74,277$  & $47,527$  & $70,028$  & $78,365$  \\
			$6$ & $>10^6(4)$ & $>10^6(4)$ & $562,156$ & $932,395$   & $260,793$ & $126,061$ & $141,914$ & $\mathbf{86,105}$  & $907,497$ & $662,983$ & $345,082$ & $280,069$ \\
			$7$ & $266,007$  & $220,647$  & $77,998$  & $317,351$   & $110,237$ & $33,691$  & $\mathbf{27,380}$  & $41,751$  & $86,269$  & $135,889$ & $54,400$  & $145,022$ \\
			$8$ & $>10^6(8)$ & $>10^6(6)$ & $776,052$ & $>10^6(18)$ & $472,125$ & $229,583$ & $313,420$ & $\mathbf{210,483}$ & $960,573$ & $700,615$ & $436,072$ & $718,803$ \\
			\bottomrule
	\end{tabular}}
	\label{tab:results4}
\end{sidewaystable}

\begin{figure}
	\resizebox{\textwidth}{!}{
		\begin{tikzpicture}
			\begin{axis}[
				legend pos=north west,
				title  = {Operational characteristics},
				xlabel = {Function evaluations},
				xmode=log,
				ymin=-0.02,ymax=1.02,
				ytick distance=0.1,
				xmode=log,
				xmin=10,
				xmax=1000001,
				xtick distance=10,
				ylabel = {Proportion of solved problems},
				ylabel style={yshift=-0.5em},
				legend style={font=\tiny,xshift=-0.5em},
				legend cell align={left},
				legend columns=1,
				height=0.75\textwidth,width=\textwidth,
				every axis plot/.append style={very thick},
				]
				\node [draw,fill=white] at (rel axis cs: 0.85,0.14) {\shortstack[l]{
						{\scriptsize \textbf{Partitioning  schemes}} \\
						\ref{p1} {\scriptsize N-DTC} \\
						\ref{p2} {\scriptsize 1-DTC} \\
						\ref{p3} {\scriptsize 1-DBDP} \\
						\ref{p4} {\scriptsize 1-DTDV}}};

				\addplot[postaction={decoration={markings,mark=between positions 0 and 1 step 0.06 with {\node[circle,draw=black,fill=princetonorange,inner sep=1.5pt,solid] {};}},decorate,},sandstorm,line width=0.75pt,densely dashed] table[x=T,y=DDA] {data/Overallassw.txt};
				\addplot[postaction={decoration={markings,mark=between positions 0 and 1 step 0.07 with {\node[mark=square,draw=black,fill=yellow,inner sep=1.5pt,solid] {};}},decorate,},sandstorm,line width=0.75pt,densely dashed] table[x=T,y=DRA] {data/Overallassw.txt};
				\addplot[postaction={decoration={markings,mark=between positions 0 and 1 step 0.08 with 	{\node[diamond,draw=black,fill=sienna,inner sep=1.5pt,solid] {};}},decorate,},sandstorm,line width=0.75pt,densely dashed] table[x=T,y=BIA] {data/Overallassw.txt};
				\addplot[postaction={decoration={markings,mark=between positions 0 and 1 step 0.09 with {\node[regular 	polygon,regular polygon sides=3,draw=black,fill=psychedelicpurple,inner sep=1pt,solid] {};}},decorate,},sandstorm,line width=0.75pt,densely dashed] table[x=T,y=ADA] {data/Overallassw.txt};
				
				\addplot[postaction={decoration={markings,mark=between positions 0 and 1 step 0.1 with {\node[circle,draw=black,fill=princetonorange,inner sep=1.5pt] {};}},decorate,},blue,line width=0.75pt] table[x=T,y=DDO] {data/Overallassw.txt};
				\addplot[postaction={decoration={markings,mark=between positions 0 and 1 step 0.11 with {\node[mark=square,draw=black,fill=yellow,inner sep=1.5pt] {};}},decorate,},blue,line width=0.75pt] table[x=T,y=DRO] {data/Overallassw.txt};
				\addplot[postaction={decoration={markings,mark=between positions 0 and 1 step 0.12 with 	{\node[diamond,draw=black,fill=sienna,inner sep=1.5pt] {};}},decorate,},blue,line width=0.75pt] table[x=T,y=BIO] {data/Overallassw.txt};
				\addplot[postaction={decoration={markings,mark=between positions 0 and 1 step 0.13 with {\node[regular 	polygon,regular polygon sides=3,draw=black,fill=psychedelicpurple,inner sep=1pt] {};}},decorate,},blue,line width=0.75pt] table[x=T,y=ADO] {data/Overallassw.txt};
				
				\addplot[postaction={decoration={markings,mark=between positions 0 and 1 step 0.14 with {\node[circle,draw=black,fill=princetonorange,inner sep=1.5pt,solid] {};}},decorate,},onyx,line width=0.75pt,densely dotted] table[x=T,y=DDG] {data/Overallassw.txt};
				\addplot[postaction={decoration={markings,mark=between positions 0 and 1 step 0.15 with {\node[mark=square,draw=black,fill=yellow,inner sep=1.5pt,solid] {};}},decorate,},onyx,line width=0.75pt,densely dotted] table[x=T,y=DRG] {data/Overallassw.txt};
				\addplot[postaction={decoration={markings,mark=between positions 0 and 1 step 0.16 with 	{\node[diamond,draw=black,fill=sienna,inner sep=1.5pt,solid] {};}},decorate,},onyx,line width=0.75pt,densely dotted] table[x=T,y=BIG] {data/Overallassw.txt};
				\addplot[postaction={decoration={markings,mark=between positions 0 and 1 step 0.17 with {\node[regular 	polygon,regular polygon sides=3,draw=black,fill=psychedelicpurple,inner sep=1pt,solid] {};}},decorate,},onyx,line width=0.75pt,densely dotted] table[x=T,y=ADG] {data/Overallassw.txt};
				
				\node [draw,fill=white] at (rel axis cs: 0.2,0.88) {\shortstack[l]{
						{\scriptsize \textbf{POH selection schemes}} \\
						\ref{p5} {\scriptsize Improved Original (IO)} \\
						\ref{p6} {\scriptsize Improved Aggressive (IA)} \\
						\ref{p7} {\scriptsize Global-Local (GL) Pareto}}};
			\end{axis}
	\end{tikzpicture}}
	\caption{Operational characteristics for all twelve \direct-type algorithmic variations on four ``hard'' GKLS-type classes}
	\label{fig:perf-l3}
\end{figure}
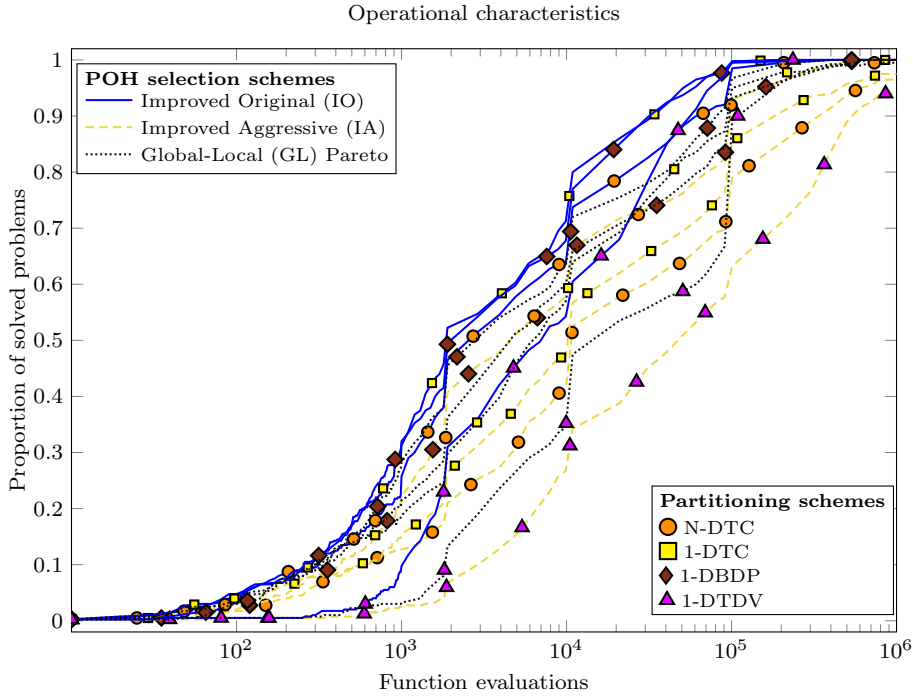

Additionally, we visualize the performance of all twelve \direct-type variations on GKLS-type problems using the operational characteristics.
\Cref{fig:perf-l2} shows the behavior on four ``simple'' GKLS classes, while \Cref{fig:perf-l3} shows the ``hard'' ones.
For simple classes, IO and GL selection schemes seem the most promising.
Among algorithms, when a low budget of function evaluations is considered (M$_{\rm max} \leq 400$), the 1-DTC-IO algorithm is the most efficient.
However, when M$_{\rm max} > 400 $, the 1-DBDP-GL algorithmic combination outperforms all others.
Furthermore, the best performance on these simple classes is achieved regardless of the selection scheme used with the 1-DBDP  partitioning strategy.

Finally, for ``hard'' classes (see \Cref{fig:perf-l3}) IO selection scheme seems the most favorable (especially for simpler problems), while GL is the second-best option.
However, when a higher maximal number of function evaluations is allowed, the performance of GL and IO selection scheme-based algorithms is quite similar.
Among the algorithms, 1-DBDP-IO and 1-DTC-IO are the two best-performing ones.

\begin{figure}
	\resizebox{\textwidth}{!}{
		\begin{tikzpicture}
			\begin{axis}[
				legend pos=north west,
				title  = {Operational characteristics},
				xlabel = {Function evaluations},
				xmode=log,
				ymin=-0.02,ymax=1.02,
				ytick distance=0.1,
				xmode=log,
				xmin=10,
				xmax=1000000,
				xtick distance=10,
				ylabel = {Proportion of solved problems},
				ylabel style={yshift=-0.5em},
				legend style={font=\tiny,xshift=-0.5em},
				legend cell align={left},
				legend columns=1,
				height=0.75\textwidth,width=\textwidth,
				every axis plot/.append style={very thick},
				]
				\node [draw,fill=white] at (rel axis cs: 0.85,0.14) {\shortstack[l]{
						{\scriptsize \textbf{Partitioning  schemes}} \\
						\ref{p1} {\scriptsize N-DTC} \\
						\ref{p2} {\scriptsize 1-DTC} \\
						\ref{p3} {\scriptsize 1-DBDP} \\
						\ref{p4} {\scriptsize 1-DTDV}}};
				
				\addplot[postaction={decoration={markings,mark=between positions 0 and 1 step 0.06 with {\node[circle,draw=black,fill=princetonorange,inner sep=1.5pt,solid] {};}},decorate,},sandstorm,line width=0.75pt,densely dashed] table[x=T,y=DDA] {data/Overallass1.txt};
				\addplot[postaction={decoration={markings,mark=between positions 0 and 1 step 0.07 with {\node[mark=square,draw=black,fill=yellow,inner sep=1.5pt,solid] {};}},decorate,},sandstorm,line width=0.75pt,densely dashed] table[x=T,y=DRA] {data/Overallass1.txt};
				\addplot[postaction={decoration={markings,mark=between positions 0 and 1 step 0.08 with 	{\node[diamond,draw=black,fill=sienna,inner sep=1.5pt,solid] {};}},decorate,},sandstorm,line width=0.75pt,densely dashed] table[x=T,y=BIA] {data/Overallass1.txt};
				\addplot[postaction={decoration={markings,mark=between positions 0 and 1 step 0.09 with {\node[regular 	polygon,regular polygon sides=3,draw=black,fill=psychedelicpurple,inner sep=1pt,solid] {};}},decorate,},sandstorm,line width=0.75pt,densely dashed] table[x=T,y=ADA] {data/Overallass1.txt};
				
				\addplot[postaction={decoration={markings,mark=between positions 0 and 1 step 0.1 with {\node[circle,draw=black,fill=princetonorange,inner sep=1.5pt] {};}},decorate,},blue,line width=0.75pt] table[x=T,y=DDO] {data/Overallass1.txt};
				\addplot[postaction={decoration={markings,mark=between positions 0 and 1 step 0.11 with {\node[mark=square,draw=black,fill=yellow,inner sep=1.5pt] {};}},decorate,},blue,line width=0.75pt] table[x=T,y=DRO] {data/Overallass1.txt};
				\addplot[postaction={decoration={markings,mark=between positions 0 and 1 step 0.12 with 	{\node[diamond,draw=black,fill=sienna,inner sep=1.5pt] {};}},decorate,},blue,line width=0.75pt] table[x=T,y=BIO] {data/Overallass1.txt};
				\addplot[postaction={decoration={markings,mark=between positions 0 and 1 step 0.13 with {\node[regular 	polygon,regular polygon sides=3,draw=black,fill=psychedelicpurple,inner sep=1pt] {};}},decorate,},blue,line width=0.75pt] table[x=T,y=ADO] {data/Overallass1.txt};
				
				\addplot[postaction={decoration={markings,mark=between positions 0 and 1 step 0.14 with {\node[circle,draw=black,fill=princetonorange,inner sep=1.5pt,solid] {};}},decorate,},onyx,line width=0.75pt,densely dotted] table[x=T,y=DDG] {data/Overallass1.txt};
				\addplot[postaction={decoration={markings,mark=between positions 0 and 1 step 0.15 with {\node[mark=square,draw=black,fill=yellow,inner sep=1.5pt,solid] {};}},decorate,},onyx,line width=0.75pt,densely dotted] table[x=T,y=DRG] {data/Overallass1.txt};
				\addplot[postaction={decoration={markings,mark=between positions 0 and 1 step 0.16 with 	{\node[diamond,draw=black,fill=sienna,inner sep=1.5pt,solid] {};}},decorate,},onyx,line width=0.75pt,densely dotted] table[x=T,y=BIG] {data/Overallass1.txt};
				\addplot[postaction={decoration={markings,mark=between positions 0 and 1 step 0.17 with {\node[regular 	polygon,regular polygon sides=3,draw=black,fill=psychedelicpurple,inner sep=1pt,solid] {};}},decorate,},onyx,line width=0.75pt,densely dotted] table[x=T,y=ADG] {data/Overallass1.txt};
				\node [draw,fill=white] at (rel axis cs: 0.2,0.88) {\shortstack[l]{
						{\scriptsize \textbf{POH selection schemes}} \\
						\ref{p5} {\scriptsize Improved Original (IO)} \\
						\ref{p6} {\scriptsize Improved Aggressive (IA)} \\
						\ref{p7} {\scriptsize Global-Local (GL) Pareto}}};
			\end{axis}
	\end{tikzpicture}}
	\caption{Operational characteristics for all twelve \direct-type algorithmic variations on four ``simple'' GKLS-type classes}
	\label{fig:perf-l2}
\end{figure}
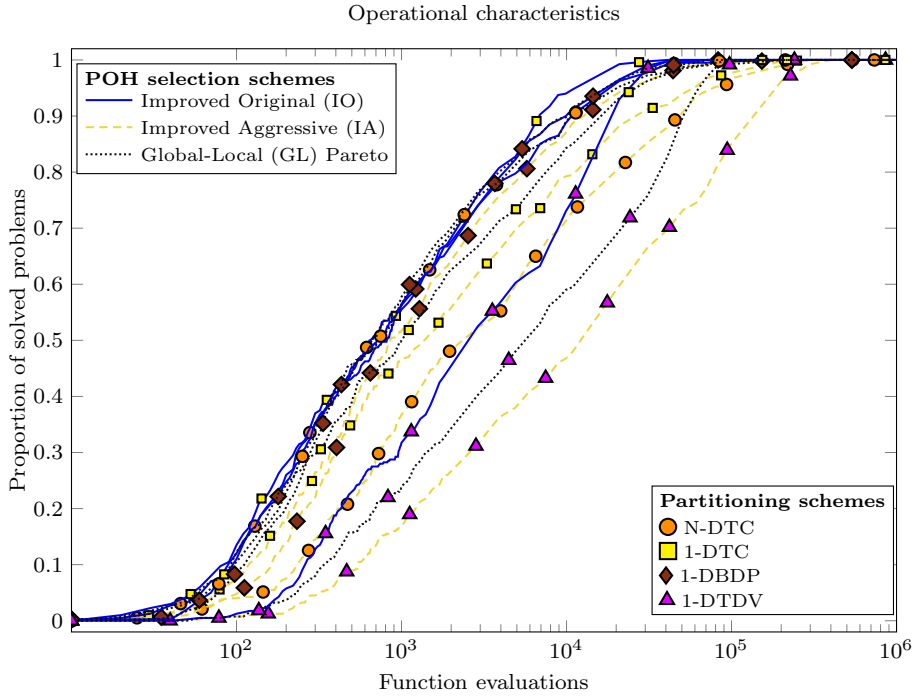

\section{Conclusions and future work}
\label{sec:conclusiuo}

This paper presented an extensive experimental investigation of various candidate selection and partitioning techniques traditionally used in the \direct-type algorithms.
Twelve \direct-type algorithmic combinations were created by considering four well-known partitioning and three selection schemes.
In general, experimental results confirmed the well-known fact from ``No free lunch theorems for optimization''~\cite{Wolpert1997} that no one particular optimization algorithm works best for every problem.
However, detailed experimental studies have helped identify particular \direct-type algorithmic variations that work well in certain situations.
For example, our experimental findings in \Cref{ssec:perturbation} revealed that what initially looks like a clear dominance case goes away with small domain perturbations.
This should remind us how dangerous it can be to generalize from limited test-function results.
Below, we emphasize when certain variations have performed best and make some recommendations based on that.

Investigation using \directlib{} test problems showed that independently on the partitioning strategy, a two-step-based Pareto selection scheme (GL) ensures the best performance on more challenging optimization problems (higher-dimensionality, multi-modal, non-convex).
The two best algorithmic variations are when the (GL) selection scheme is combined with the 1-DTC and 1-DBDP partitioning approaches.
While the 1-DTDV-GL looks best for more straightforward problems (low-dimensional, uni-modal), the 1-DTC-GL, 1-DBDP-GL combination is more efficient in solving more challenging problems.
The worst results were obtained using various partitioning strategies combined with the (IO) selection scheme, which showed promising performance only when a given budget of function evaluations is small.

Moreover, regardless of the selection scheme, the 1-DTDV partitioning strategy has a signiﬁcant advantage when the most solution coordinates are on the boundary of the feasible region.
Additionally, the 1-DTDV partitioning approach has proven to be the most efficient in solving low-dimensional \directlib{} test problems.
However, the combination based on the 1-DTDV partitioning scheme is very inefficient for higher dimensional test problems.
For such problems, the 1-DBDP partitioning approach seems much more appropriate.


Experimental investigation on 800 GKLS-type test problems showed contrasting results.
This study revealed that the (IO) scheme could be very efficient.
While on simple GKLS classes, the efficiency of \direct-type algorithms based on (IO) and (GL) selection schemes are very similar.
However, better performance is explicitly achieved with (IO) for hard classes.
Let us recall that this selection scheme showed the worst results in the previous investigation.

To sum up, our study demonstrated that using already known techniques combined in new variations can create more efficient \direct-type algorithms.
For example, efficient diagonal partition-based \birect{} can be further improved by replacing the original selection scheme with the GL selection (from the \directgl{} algorithm), resulting in a more efficient algorithm (1-DBDP-GL).


As for further research, one possible direction could be improving the two-step-based (Global-Local) Pareto selection scheme (GL). Algorithms based on this scheme showed superior performance solving most complex optimization test problems but relatively poor efficiency on more straightforward problems.
One possible modification could be borrowing \direct's $\varepsilon$ parameter or similar technique to limit the size of selected POHs.
Optionally, instead of performing the selection enhancing the local search in every iteration, a specific rule could be added to specify when this selection is needed.

Finally, finding the solution efficiently should start by investigating the problem.
Then, considering this knowledge, the design or finding of a specific optimization algorithm is needed.
Thus, one of our nearest future work plans is to extend this idea by developing the automatic \direct-type algorithm selection.

\section*{Source code statement}
All twelve introduced \direct-type algorithms are implemented in \texttt{MATLAB} and are available in the most recent version of \directgo{} (\url{https://github.com/blockchain-group/DIRECTGO/tree/v1.1.0}) and can be used under the MIT license.
We welcome contributions and corrections to it.

\section*{Data statement}
\texttt{DIRECTGOLib} - \direct{} \textbf{G}lobal \textbf{O}ptimization test problems \textbf{Lib}rary is designed as a continuously-growing open-source GitHub repository to which anyone can easily contribute.
The exact data underlying this article from \directlib{} can be accessed either on GitHub or at Zenodo:
\begin{itemize}
  \item at GitHub: \url{https://github.com/blockchain-group/DIRECTGOLib/tree/v1.1},
  \item at Zenodo: \url{https://doi.org/10.5281/zenodo.6491951},
\end{itemize}
and used under the MIT license.
We welcome contributions and corrections to it.

\vspace{15pt}
\noindent\textbf{Funding} The research work of S. Stripinis was funded by a Grant (No. S-MIP-21-53) from the \textit{Research Council of Lithuania}.

\vspace{15pt}
\noindent\textbf{Acknowledgment} The authors greatly thank the anonymous Reviewer for his valuable and constructive comments, which helped us significantly extend and improve the manuscript.

\bibliographystyle{spmpsci}      
\bibliography{library}   

\appendix

\section{\directlib{} library}
\label{apendixas}

A summary of all used box-constrained optimization problems from \directlib~\cite{DIRECTGOLib2022v11} and their properties are given in \Cref{tab:test}. 
Here, the main features are reported: problem number (\#), name of the problem, source, dimensionality ($ n $), default optimization domain ($ D $), perturbed optimization domain ($\tilde{D}$), problem type, and the known minimum ($ f^* $).
The default domains are taken from the literature.
Whenever the global minimum point lies at the initial sampling points for at least one tested algorithm, the domain $D$ has been shifted by $22.5 \%$ to the right side.
These modified problems are marked with the beta sign $(^{\beta})$.
Here the sign ``-'' means that $\tilde{D}$ is the same as $D$.
Some of these test problems have several variants, e.g., \textit{Bohachevsky}, \textit{Shekel}, and some of them, like \textit{Alpine}, \textit{Csendes}, \textit{Griewank}, etc., can be tested for varying dimensionality.

\begin{table}
	\caption{Key characteristics of the \directlib~\cite{DIRECTGOLib2022v11} test problems for box-constrained global optimization}
	\resizebox{\textwidth}{!}{
		\begin{tabular}[tb]{@{\extracolsep{\fill}}rlrrrrrrr}
			\toprule
			\# & Name & Source & $n$  &  $ D $ &  $ \tilde{D} $ & Type & No. of minima & $ f^* $ \\
			\midrule
			$1,2,3$    & \textit{Ackley}$^{\beta}$ 		& \cite{Hedar2005,Derek2013} 		& $2,5,10$ 	& $[-15, 35]^n$ & $[-18, 47]^n$ 			    		& non-convex & multi-modal & $0.0000$			\\
			$4,5,6$    & \textit{Alpine}$^{\alpha}$ 		& \cite{Gavana2021} 		   		& $2,5,10$ 	& $[0, 10]^n$   & $[\sqrt[i]{2}, 8 + \sqrt[i]{2}]^n$  	& non-convex & multi-modal & $-2.8081^n$ 		\\
			$7$ 	   & \textit{Beale}      	 			& \cite{Hedar2005,Derek2013} 		& $2$  		& $[-4.5, 4.5]^n$ & $-$  	   		 		& non-convex & multi-modal & $0.0000$   		\\
			$8$ 	   & \textit{Bohachevsky$1^{\beta}$}  	& \cite{Hedar2005,Derek2013} 		& $2$  		& $[-100, 110]^n$ & $[-55, 145]^n$   	   							& convex 	 & uni-modal   & $0.0000$   		\\
			$9$ 	   & \textit{Bohachevsky$2^{\beta}$} 	& \cite{Hedar2005,Derek2013} 		& $2$  		& $[-100, 110]^n$ & $[-55, 145]^n$   	   							& non-convex & multi-modal & $0.0000$   		\\
			$10$ 	   & \textit{Bohachevsky$3^{\beta}$}  	& \cite{Hedar2005,Derek2013} 		& $2$  		& $[-100, 110]^n$ & $[-55, 145]^n$   	   							& non-convex & multi-modal & $0.0000$   		\\
			$11$ 	   & \textit{Booth}          			& \cite{Hedar2005,Derek2013} 		& $2$  		& $[-10, 10]^n$ & $-$  	 	 						& convex 	 & uni-modal   & $0.0000$			\\
			$12$ 	   & \textit{Branin}          			& \cite{Hedar2005,Dixon1978} 		& $2$  		& $[-5, 10] \times [10,15]$ & $-$ 					& non-convex & multi-modal & $0.3978$			\\
			$13$ 	   & \textit{Bukin6}           			& \cite{Derek2013} 		  	 		& $2$  		& $[-15, 5] \times [-3,3]$ & $-$  					& convex 	 & multi-modal & $0.0000$   		\\
			$14$ 	   & \textit{Colville}        			& \cite{Hedar2005,Derek2013} 		& $4$  		& $[-10, 10]^n$ & $-$           						& non-convex & multi-modal & $0.0000$   		\\
			$15$ 	   & \textit{Cross\_in\_Tray}   		& \cite{Derek2013} 		  	 		& $2$    	& $[0, 10]^n$ & $-$           						& non-convex & multi-modal & $-2.0626$   		\\
			$16$ 	   & \textit{Crosslegtable}   			& \cite{Gavana2021} 		  	 	& $2$    	& $[-10, 15]^n$ & $-$           						& non-convex & multi-modal & $-1.000$   		\\
			$17,18,19$ & \textit{Csendes}$^{\beta}$     	& \cite{Gavana2021} 		  	 	& $2,5,10$  & $[-10, 21]^n$ & $[-10, 25]^n$          	& convex 	 & multi-modal & $0.0000$   		\\
			$20$ 	   & \textit{Damavandi}   				& \cite{Gavana2021} 		  	 	& $2$    	& $[0, 14]^n$  & $-$           						& non-convex & multi-modal & $0.0000$   		\\
			$21,22,23$ & \textit{Deb$01$}$^{\beta}$    	& \cite{Gavana2021} 				& $2,5,10$  & $[-1, 1]^n$ & $[-0.55, 1.45]^n$         					& non-convex & multi-modal & $-1.0000$   		\\
			$24,25,26$ & \textit{Deb$02$}$^{\beta}$   		& \cite{Gavana2021} 				& $2,5,10$  & $[0, 1]^n$ & $[0.225, 1.225]^n$ 							& non-convex & multi-modal & $-1.0000$   		\\
			$27,28,29$ & \textit{Dixon\_and\_Price}   		& \cite{Hedar2005,Derek2013} 		& $2,5,10$  & $[-10, 10]^n$ & $-$         						& convex 	 & multi-modal & $0.0000$   		\\
			$30$ 	   & \textit{Drop\_wave}$^{\beta}$   	& \cite{Derek2013} 		  	 		& $2$    	& $[-5.12, 6.12]^n$  & $[-4, 6]^n$      								& non-convex & multi-modal & $-1.0000$   		\\
			$31$ 	   & \textit{Easom}$^{\alpha}$     		& \cite{Hedar2005,Derek2013} 		& $2$    	& $[-100, 100]^n$ & $\left[\dfrac{-100}{i+1}, 100i \right]^n$ 	& non-convex & multi-modal & $-1.0000$			\\
			$32$ 	   & \textit{Eggholder}       			& \cite{Derek2013} 		  	 		& $2$   	& $[-512, 512]^n$ & $-$         						& non-convex & multi-modal & $-959.6406$	 	\\
			$33$ 	   & \textit{Goldstein\_and\_Price}$^{\beta}$ & \cite{Hedar2005,Dixon1978} & $2$   	& $[-2, 2]^n$   & $[-1.1, 2.9]^n$          						& non-convex & multi-modal & $3.0000$ 			\\
			$34,35,36$ & \textit{Griewank}$^{\alpha}$     & \cite{Hedar2005,Derek2013} & $2,5,10$ & $[-600, 700]^n$ & $\left[-\sqrt{600i}, \dfrac{600}{\sqrt{i}} \right]^n$ & non-convex & multi-modal & $0.0000$   		\\
			$37$ 	   & \textit{Hartman$3$}       			& \cite{Hedar2005,Derek2013} 		& $3$    	& $[0, 1]^n$ 	& $-$              						& non-convex & multi-modal & $-3.8627$			\\
			$38$ 	   & \textit{Hartman$6$}        		& \cite{Hedar2005,Derek2013} 		& $6$    	& $[0, 1]^n$  & $-$            						& non-convex & multi-modal & $-3.3223$			\\
			$39$ 	   & \textit{Holder\_Table}      		& \cite{Derek2013} 			 		& $2$    	& $[-10, 10]^n$   & $-$           						& non-convex & multi-modal & $-19.2085$			\\
			$40$ 	   & \textit{Hump}              		& \cite{Hedar2005,Derek2013} 		& $2$    	& $[-5, 5]^n$ & $-$           						& non-convex & multi-modal & $-1.0316$			\\
			$41$ 	   & \textit{Langermann}        		& \cite{Derek2013} 			 		& $2$    	& $[0, 10]^n$   & $-$           						& non-convex & multi-modal & $-4.1558$			\\
			$42,43,44$ & \textit{Levy}              		& \cite{Hedar2005,Derek2013} 		& $2,5,10$  & $[-5, 5]^n$ & $-$             						& non-convex & multi-modal & $0.0000$   		\\
			$45$ 	   & \textit{Matyas}$^{\beta}$         & \cite{Hedar2005,Derek2013} 		& $2$    	& $[-10, 15]^n$ & $[-5.5, 14.5]^n$  	 						& convex 	 & uni-modal   & $0.0000$ 			\\
			$46$ 	   & \textit{McCormick}         		& \cite{Derek2013} 	 		 		& $2$    	& $[-1.5, 4] \times [-3,4]$ & $-$ 					& convex 	 & multi-modal & $-1.9132$ 			\\
			$47$ 	   & \textit{Michalewicz}       		& \cite{Hedar2005,Derek2013} 		& $2$    	& $[0, \pi]^n$ & $-$            						& non-convex & multi-modal & $-1.8013$	 		\\
			$48$ 	   & \textit{Michalewicz}       		& \cite{Hedar2005,Derek2013} 		& $5$    	& $[0, \pi]^n$ & $-$          						& non-convex & multi-modal & $-4.6876$			\\
			$49$     & \textit{Michalewicz}       		& \cite{Hedar2005,Derek2013} 		& $10$   	& $[0, \pi]^n$ & $-$         	 						& non-convex & multi-modal & $-9.6601$			\\
			$50$ 	   & \textit{Permdb$4$} 						& \cite{Hedar2005,Derek2013} 	& $4$    	& $[-i, i]^n$ & $[-i, i]^n$  									& non-convex & multi-modal & $0.0000$   		\\
			$51,52,53$ & \textit{Pinter}$^{\beta}$     	& \cite{Gavana2021} 		  	 	& $2,5,10$  & $[-10, 10]^n$ & $[-5.5, 14.5]^n$      						& non-convex & multi-modal & $0.0000$   		\\
			$54$ 	   & \textit{Powell}       				& \cite{Hedar2005,Derek2013} 		& $4$    	& $[-4, 5]^n$ & $-$             						& convex 	 & multi-modal & $0.0000$   		\\
			$55$ 	   & \textit{Power\_Sum}$^{\alpha}$     & \cite{Hedar2005,Derek2013} 		& $4$    	& $[0, 4]^n$ & $[1, 4 + \sqrt[i]{2}]^n$            			& convex 	 & multi-modal   & $0.0000$   		\\
			$56,57,58$ & \textit{Qing}     	    			& \cite{Gavana2021} 		  	 	& $2,5,10$  & $[-500, 500]^n$ & $-$      	 						& non-convex & multi-modal & $0.0000$   		\\
			$59,60,61$ & \textit{Rastrigin}$^{\alpha}$    	& \cite{Hedar2005,Derek2013} 		& $2,5,10$  & $[-6.12, 5.12]^n$ & $[-5\sqrt[i]{2}, 7+\sqrt[i]{2}]^n$       		& non-convex & multi-modal & $0.0000$   		\\
			$62,63,64$ & \textit{Rosenbrock}$^{\alpha}$   	& \cite{Hedar2005,Dixon1978}  		& $2,5,10$  & $[-5, 10]^n$ & $\left[-\dfrac{5}{\sqrt{i}}, 10\sqrt{i} \right]^n$ & non-convex & uni-modal & $0.0000$   		\\
			$65,66,67$ & \textit{Rotated\_H\_Ellip}$^{\beta}$ & \cite{Derek2013} 	 			& $2,5,10$  & $[-65.536, 66.536]^n$ & $[-35, 96]^n$   								& convex 	 & uni-modal   & $0.0000$   		\\
			$68,69,70$ & \textit{Schwefel}$^{\alpha}$ 		& \cite{Hedar2005,Derek2013} 		& $2,5,10$ 	& $[-500, 500]^n$ 	& $\left[-500 + \dfrac{100}{\sqrt{i}}, 500 - \dfrac{40}{\sqrt{i}} \right]^n$ & non-convex & multi-modal & $0.0000$ \\
			$71$ 	   & \textit{Shekel$5$}  				& \cite{Hedar2005,Derek2013} 		& $4$    	& $[0, 10]^n$ 		& $-$             						& non-convex & multi-modal & $-10.1531$			\\
			$72$ 	   & \textit{Shekel$7$}   				& \cite{Hedar2005,Derek2013}	 	& $4$    	& $[0, 10]^n$ 		& $-$             						& non-convex & multi-modal & $-10.4029$			\\
			$73$ 	   & \textit{Shekel$10$}  				& \cite{Hedar2005,Derek2013} 		& $4$    	& $[0, 10]^n$ 		& $-$             						& non-convex & multi-modal & $-10.5364$	 		\\
			$74$ 	   & \textit{Shubert}           		& \cite{Hedar2005,Derek2013} 		& $2$    	& $[-10, 10]^n$ 	& $-$           						& non-convex & multi-modal & $-186.7309$		\\
			$75,76,77$ & \textit{Sphere}$^{\beta}$      	& \cite{Hedar2005,Derek2013} 		& $2,5,10$  & $[-5.12, 6.12]^n$ & $[-2.75, 7.25]^n$      						& convex 	 & uni-modal   & $0.0000$   		\\
			$78,79,80$ & \textit{Styblinski\_Tang}$^{\alpha}$ & \cite{Clerc1999} 			    & $2,5,10$  & $[-5, 5]^n$ 		& $[-5, 5+\sqrt[i]{3}]^n$       	    					& non-convex & multi-modal & $-39.1661n$		\\
			$81,82,83$ & \textit{Sum\_of\_Powers}$^{\beta}$ & \cite{Derek2013}		 			& $2,5,10$  & $[-1, 2.5]^n$ 		& $[-0.55, 1.45]^n$      						& convex 	 & uni-modal   & $0.0000$   		\\
			$84,85,86$ & \textit{Sum\_Square}$^{\beta}$  	& \cite{Clerc1999} 			 		& $2,5,10$  & $[-10, 15]^n$ 	& $[-5.5, 14.5]^n$        						& convex 	 & uni-modal   & $0.0000$			\\
			$87$ 	   & \textit{Trefethen}           		& \cite{Gavana2021} 				& $2$    	& $[-2, 2]^n$ 		& $-$           						& non-convex & multi-modal & $-3.3068$			\\
			$88,89,90$ & \textit{Trid}              		& \cite{Hedar2005,Derek2013} 		& $2,5,10$  & $[-100, 100]^n$  	& $-$           					& convex 	 & multi-modal & $\vartheta$		\\
			$91,92,93$ & \textit{Vincent}           		& \cite{Clerc1999} 					& $2,5,10$  & $[0.25, 10]^n$ 	& $-$          						& non-convex & multi-modal & $-n$   			\\
			$94,95,96$ & \textit{Zakharov}$^{\beta}$     	& \cite{Hedar2005,Derek2013} 		& $2,5,10$  & $[-5, 11]^n$ 		& $[-1.625, 13.375]^n$      					& convex 	 & multi-modal & $0.0000$    		\\
			\midrule
			\multicolumn{8}{l}{$\vartheta$ -- $-\frac{1}{6}n^3 - \frac{1}{2}n^2 +  \frac{2}{3}n$}\\
			\multicolumn{8}{l}{${\alpha}$ -- domain $D$ was perturbed to avoid the dominance of certain partitioning schemes (see \Cref{ssec:perturbation})} \\
			\multicolumn{8}{l}{${\beta}$ -- domain $D$ was moved by $22.5 \%$ percent to the right side} \\
			\multicolumn{8}{l}{${i = 1, ..., n}$} \\
			
			\bottomrule
	\end{tabular}}
	\label{tab:test}
\end{table}

\end{document}